\documentclass[a4paper]{amsart}
\usepackage{mathtools,amsfonts,amssymb,mathrsfs}
\usepackage[a4paper,hmargin=26mm,vmargin=28mm]{geometry}
\usepackage{bm,multicol,mathbbol}
\usepackage{booktabs}
\usepackage{anyfontsize}
\usepackage{calligra}
\usepackage[svgnames]{xcolor}
\usepackage[pagebackref=true]{hyperref}
\usepackage{tikz}
\usetikzlibrary{matrix}
\usepackage{enumitem}
\setlist[enumerate]{nosep}
\usepackage[hyperpageref]{backref}

\title{Fayers' conjecture and the socles of cyclotomic Weyl modules}
\subjclass[2000]{20G43, 20C08, 20C30, 05E10}
\keywords{Cyclotomic Hecke algebras, Schur algebras, quasi-hereditary and graded
cellular algebras, Khovanov--Lauda--Rouquier algebras}
\author{Jun Hu}\address{School of Mathematics,
Beijing Institute of Technology,
Beijing, 100081, P.R.~China}

\def\Email#1{\email{\href{mailto:#1}{#1}}}
\Email{junhu303@qq.com}

\author{Andrew Mathas}
\address{School of Mathematics and Statistics, University of Sydney, NSW 2006, Australia}
\Email{andrew.mathas@sydney.edu.au}
\hypersetup{
  pdftitle    ={Fayers conjecture and the socles of cyclotomic Weyl modules},%
  pdfauthor   ={Jun Hu and Andrew Mathas},%
  colorlinks  = true,%
  linkcolor   = blue,%
  anchorcolor = red,%
  citecolor   = blue,%
  urlcolor    = blue%
}

\renewcommand*{\backref}[1]{}
\renewcommand*{\backrefalt}[4]{{[\tiny%
  \ifcase #1 Not cited.\relax\or Page~#2.%
  \else Pages #2.\fi]}
}

\DeclareFontFamily{OT1}{pzc}{}
\DeclareFontShape{OT1}{pzc}{m}{it}{<-> s * [1.20] pzcmi7t}{}
\DeclareMathAlphabet{\mathpzc}{OT1}{pzc}{m}{it}

\DeclareSymbolFontAlphabet{\mathbb}{AMSb}

\usepackage{cite}

\newcommand\Comment[2][Hu]{\space\par\medskip\noindent%
   \fbox{\begin{minipage}{\textwidth}\textbf{Comment\ifx\relax#1\else---#1\fi}\newline%
        #2\end{minipage}}\medskip
}
\usepackage{tikz}
\usetikzlibrary{matrix}

\let\<=\langle
\let\>=\rangle
\def\({\big(}
\def\){\big)}

\newcommand\bijection[1][\sim]{\overset{#1}{\longrightarrow}}
\let\surjection\twoheadrightarrow
\let\injection\hookrightarrow

\newcommand{\C}{\mathbb{C}}
\newcommand{\N}{\mathbb N}

\newcommand{\Q}{\mathbb Q}
\newcommand{\Z}{\mathbb Z}
\newcommand{\Sym}{\mathfrak S}

\def\map#1#2{\,{:}\,#1\!\longrightarrow\!#2}

\usepackage{xparse,mathtools}
\DeclarePairedDelimiterX{\set}[1]{\{}{\}}{\setargs{#1}}
\NewDocumentCommand{\setargs}{>{\SplitArgument{1}{|}}m}{\setargsaux#1}
\NewDocumentCommand{\setargsaux}{mm}
{\IfNoValueTF{#2}{#1} {#1\,\delimsize|\,\mathopen{}#2}}

\def\pmod#1{\text{ }(\text{mod } #1)\,}

\def\K{\mathscr{K}}
\def\P{\mathscr{P}}

\def\Fun{\mathsf{F}}
\newcommand\Funn[1][n]{\Fun_{#1}}
\newcommand\HEquiv[1][n]{\mathsf{E}^\#_{#1}}
\newcommand\Run[1][n]{\mathsf{R}_{#1}}

\newcommand\HH{\mathcal{H}}   
\newcommand\Hn[1][n]{\HH_{#1}}

\newcommand\Hlam[1][\blam]{\Hn^{\gdom#1}}
\newcommand\Hllam[1][\blam]{\Hn^{\ldom#1}}

\newcommand\R[1][n]{\RR_{#1}}

\newcommand\RingelDual[1]{#1^{\text{RD}}}

\renewcommand\sl[1][e]{\widehat{\mathfrak{sl}}_{#1}}

\newcommand\Usl[1][e]{\mathcal{U}_q(\sl[#1])}

\newcommand\Fock[1][\charge']{\mathcal{F}({#1})}
\newcommand\BarInvolution{\rule[1.2ex]{.6em}{.1ex}}

\newcommand\DD{\mathbb{D}_\prime}
\newcommand\DL{\mathbb{L}_\prime}
\newcommand\DS{\mathbb{S}_\prime}
\newcommand\DP{\mathbb{P}_\prime}
\newcommand\DDelta{%
  \mathbb{W}_\prime
}

\def\SS{\mathcal{S}}  
\def\Sn{\SS_n}
\newcommand\Slam[1][\blam]{\SS_n^{\gdom#1}}

\newcommand\Mr{M^{\oplus r}}

\newcommand\Emunu[1][\bnu]{\mathscr{E}_{\bmu#1}}
\newcommand\ELam{E^\Lambda}

\newcommand\DSn{\dot\SS_n'}
\newcommand\DHn{\dot\HH_n'}

\def\P{\mathscr{P}}
\newcommand{\Klesh}[1][n]{\mathcal{K}_{#1}}
\newcommand{\Parts}[1][n]{\P_{#1}}
\newcommand\blam{{\boldsymbol\lambda}}
\newcommand\brho{{\boldsymbol\rho}}
\newcommand\bmu{{\boldsymbol\mu}}
\newcommand\bnu{{\boldsymbol\nu}}
\newcommand\bal{{\boldsymbol\alpha}}

\newcommand\bsig{{\boldsymbol\sigma}}
\newcommand\btau{{\boldsymbol\tau}}
\newcommand\bom{{\boldsymbol\omega}}
\newcommand\lam{{\lambda}}
\newcommand\bzero{\mathbf{\underline{0}}}

\newcommand\Dlam[1][\blam]{\Delta^{#1}}
\newcommand\ilam[1][\blam]{\iota_{#1}}
\newcommand\pilam[1][\blam]{\pi_{#1}}
\newcommand\pimu{\pilam[\bmu]}
\newcommand\plam[1][\blam]{1_{#1}}
\newcommand\pmu{\plam[\bmu]}
\newcommand\pom{\plam[\bom]}
\newcommand\zlam[1][\blam]{\zeta_{#1}}
\newcommand\zmu{\zlam[\bmu]}

\newcommand\bi{\mathbf{i}}

\def\s{{\mathfrak s}}
\def\t{{\mathfrak t}}
\def\u{{\mathfrak u}}
\def\v{{\mathfrak v}}
\def\bS{{\mathtt S}}
\def\ubs{\bS^\circ}
\def\bT{{\mathtt T}}
\def\bU{{\mathtt U}}
\def\bV{{\mathtt V}}
\def\tlam{\t^\blam}
\def\tllam{\t_\blam}

\def\Tlam{\bT^\blam}
\def\tmu{\t^\bmu}
\newcommand\phiUV[1][\bU\bV]{\varphi_{#1}}
\newcommand\phiST[1][\bS\bT]{\varphi_{#1}}

\def\Std{\mathop{\rm Std}\nolimits}
\def\RSStd{\mathcal{T}_{\text{row}}}
\def\CSStd{\mathcal{T}_{\text{col}}}

\let\gedom=\trianglerighteq
\let\gdom=\vartriangleright

\DeclareMathOperator\Gedom{\,{\underline{\kern-.1ex{\blacktriangleright}\kern-0.1ex}}\,}

\let\ldom=\vartriangleleft

\newcommand{\charge}{{\boldsymbol{\kappa}}}

\usepackage{etoolbox}
\newcommand{\DeclareMyOperator}[1]{%
    \expandafter\DeclareMathOperator\csname #1\endcsname{\mathrm{#1}}
}
\forcsvlist{\DeclareMyOperator}{%
  End,Ext,Hom,comp,Row,Col,soc,hd,id,wt,im,Shape,rad,op,res,Add,Rem,Gr,row,col,
}

\def\L{\mathop{\mathsf{L}}\nolimits}
\def\R{\mathop{\mathsf{R}}\nolimits}
\def\LR{\mathop{\mathsf{LR}}\nolimits}
\def\RL{\mathop{\mathsf{RL}}\nolimits}

\DeclareMathOperator\Mod{{-}mod}

\newcommand\rDelMod[1]{\text{mod-}\!\mathop{\mathcal{F}_W(#1)}}
\newcommand\rNabMod[1]{\text{mod-}\!\mathop{\mathcal{F}_V(#1)}}

\DeclareMathOperator\RMod{mod{-}\!}
\DeclareMathOperator\GRMod{grmod{-}\!}

\newcommand\defect[1][]{\mathop{\rm def}\nolimits_{#1}}
\def\bm{\mathbf{m}}

\def\op{{\text{op}}}

\usepackage{aliascnt}

\def\NewTheorem#1{%
    \newaliascnt{#1}{equation}
    \newtheorem{#1}[#1]{#1}
    \aliascntresetthe{#1}
    \expandafter\def\csname #1autorefname\endcsname{#1}
}

\def\equationautorefname~#1\null{(#1)\null}

\newcounter{main}
\theoremstyle{plain}

\swapnumbers
\numberwithin{equation}{section}
\NewTheorem{Assumption}
\NewTheorem{Proposition}
\NewTheorem{Theorem}
\NewTheorem{Corollary}
\NewTheorem{Conjecture}
\NewTheorem{Lemma}
\NewTheorem{Definition}
\theoremstyle{definition}
\theoremstyle{remark}
\NewTheorem{Remark}

\theoremstyle{definition}
\NewTheorem{Example}
\AtEndEnvironment{Example}{\null\hfill$\Diamond$}%

\begin{document}

\bibliographystyle{andrew}

\begin{abstract}
  Gordon James proved that the socle of a Weyl module of a classical Schur
  algebra is a sum of simple modules labelled by $p$-restricted
  partitions. We prove an analogue of this result in the very general
  setting of ``Schur pairs''. As an application we show that the socle
  of a Weyl module of a cyclotomic $q$-Schur algebra is a sum of simple
  modules labelled by Kleshchev multipartitions and we use this result
  to prove a conjecture of Fayers that leads to an efficient LLT
  algorithm for the higher level cyclotomic Hecke algebras of type~$A$.
  Finally, we prove a cyclotomic analogue of the Carter-Lusztig theorem.
\end{abstract}

\maketitle
\setcounter{tocdepth}{1}\tableofcontents

\section{Introduction} In their landmark paper~\cite{LLT}, Lascoux,
  Leclerc and Thibon conjectured that the decomposition matrices of the
  Iwahori-Hecke algebras of type~$A$ could be computed using the
  canonical bases of the \textit{level} one Fock spaces.
  Ariki~\cite{Ariki:can} generalised, and proved, this conjecture for
  the cyclotomic Hecke algebras of arbitrary level in type~$A$. Unlike
  in level one, the canonical bases of higher level Fock spaces are
  quite difficult to compute. Building on~\cite{LLT},
  Fayers~\cite{Fayers:LLT} gave a more efficient algorithm for computing
  the canonical bases of the higher level Fock spaces. As a result of
  his calculations he made a natural conjecture for the maximal degrees
  of the LLT-polynomials, which are certain parabolic Kazhdan-Lusztig
  polynomials that are computed by Fayers' algorithm. This conjecture, if true,
  further improves the efficiency of Fayers' algorithm. This paper
  started as a project to prove Fayers' conjecture but the machinery
  that we develop to do this has several other applications.

  Fayers' conjecture is a purely combinatorial statement that gives an
  upper bound on the degrees of certain polynomials. As we will show, the
  representation theoretic significance of his conjecture gives
  information about the socle of the Weyl modules. To
  prove Fayers' conjecture, we take advantage of recent advances by
  Stroppel-Webster~\cite{StroppelWebster:QuiverSchur},
  Maksimau~\cite{Maksimau:QuiverSchur},
  Brundan-Kleshchev~\cite{BK:GradedDecomp,BK:GradedKL} and the
  authors~\cite{HuMathas:QuiverSchurI} in the graded representation
  theory of the cyclotomic Hecke and Schur algebras to connect the
  higher level LLT-polynomials with the graded decomposition numbers of
  these algebras. This machinery allows us to reduce Fayers' conjecture
  to understanding the socle of a Weyl module.

  This paper starts by analysing the endomorphism algebra $S=\End_A(M)$
  of a \textit{Schur pair} $(A, M)$, where $A$ is a self-injective
  algebra and $M$ is a faithful $A$-module. In this very general
  setting, which includes the cyclotomic Hecke algebras~$\Hn$ and the
  Schur algebras~$\Sn$ as special cases, we classify the
  projective-injective $S$-modules and show that the modules that can
  appear in the socle of the projective $S$-modules are labelled by the
  simple $A$-modules. We then proceed to apply this general theory to
  the cyclotomic Hecke algebras and Schur algebras.  Our main results
  are:
  \begin{itemize}
    \item We prove that the simple modules appearing in the socle of a
    Weyl modules are indexed by \textit{Kleshchev multipartitions}
    (\autoref{T:DeltaSocle}).
    \item We prove Fayers Conjecture (\autoref{T:Fayers}).
    \item We give a cyclotomic generalisation of the Carter-Lusztig
    theorem~\cite{CarterLusztig} (\autoref{T:CarterLusztig}).
  \end{itemize}
  Along the way we prove a number of useful new results about the
  cyclotomic Hecke and Schur algebras, such as a new description of the
  simple $\Hn$-modules (\autoref{L:SimpleSubmodules}) and an interpretation
  of Ringel duality at the level of Hecke algebras (\autoref{T:RingelDuality}).

  All of the results in this paper apply to what are traditionally
  called the degenerate and non-degenerate cyclotomic Hecke algebras and
  Schur algebras (but not to the Ariki-Koike algebras with $\xi^2=1$).
  This is because, following~\cite{HuMathas:SeminormalQuiver}, we use a
  slightly different presentation (\autoref{D:Hecke}) of the
  cyclotomic Hecke algebras that simultaneously captures the degenerate
  and non-degenerate cyclotomic Hecke algebras.

\subsection*{Acknowledgements}

  We thank the referee for their extensive comments and corrections.
  Both authors were supported by the Australian Research Council. The
  first author was also supported by the National Natural Science
  Foundation of China (No. 11525102).

\section{Schur pairs and the socle of tensor space}\label{S:SchurPairs}

  Partly inspired by the extensive literature of Schur algebras acting
  on \textit{tensor space}, such
  as~\cite{Green,DJ:qWeyl,MazorchukStroppel:SerreFunctors}, this section
  studies the endomorphism algebra of a finite dimensional faithful
  module of a finite dimensional self-injective algebra. Our main aim is
  to understand, in this general setting, what simple modules can
  appear in the socle of ``tensor space''.

  Throughout this paper we fix a field~$K$. If $A$ is a $K$-algebra let
  $A\Mod$ and $\RMod A$ be the categories of finite dimensional left and
  right $A$-modules, respectively.

  Recall that a \textbf{trace form} on an algebra $A$ is a linear map
  $\tau\map AK$
  such that $\tau(ab) = \tau(ba)$, for all $a,b\in A$. The algebra $A$ is
  \textbf{symmetric} if it has a non-degenerate symmetric bilinear form
  $\< ,\ \>\map{A\times A}K$ such that
  \[
      \<xa,y\>=\<x,ay\>,\qquad\text{for all }a,x,y\in A
  \]
  The algebra $A$ is \textbf{self-injective} if the regular
  representation of~$A$ is injective as an $A$-module. As is
  well-known~\cite[\S60]{C&R}, symmetric algebras are self-injective.

  \begin{Definition}\label{D:SchurPair}
    A \textbf{Schur pair} is an ordered pair $(A, M)$, where $A$ is a
    self-injective finite dimensional $K$-algebra and~$M$ is a faithful
    finite dimensional right $A$-module. If $(A, M)$ is a Schur pair let
    \[
        S=S_A(M)=\End_A(M)
    \]
    be the algebra of $A$-module endomorphisms of~$M$.
  \end{Definition}

  By extending $K$, if necessary, we will always assume that $A$ and $S$
  are both split over~$K$.

  An explicit example to keep in mind is the case where $A$ is the
  Iwahori-Hecke algebra of the symmetric group and $M$ is tensor
  space~\cite[\S1 and \S2]{DJ:qWeyl}, so that  $S$ is the $q$-Schur
  algebra. More generally, as we show in \autoref{P:SchurPairs}, we can
  take~$A$ to be a cyclotomic Hecke algebra of type~$A$, in which case
  there are two natural choices for~$M$.

  For the remainder of this section we assume that $(A, M)$ is a Schur
  pair and that $S=S_A(M)$. Our aim is to understand the interplay
  between the representation theory of the algebras~$A$ and~$S$. Our
  starting point is the following easy fact.

  \begin{Proposition}\label{P:Equivalence}
    Suppose that $(A, M)$ is a Schur pair.
    \begin{enumerate}
      \item There exists an integer $r\ge1$ and an $A$-module~$N$
      such that $M^{\oplus r}\cong A\oplus N$ as $A$-modules.
      \item If $r\ge1$ then $(A,M^{\oplus r})$ is a Schur pair and the
      algebras $S_A(M)$ and $S_A(M^{\oplus r})$ are Morita equivalent.
    \end{enumerate}
  \end{Proposition}

  \begin{proof}
    First consider (a). Since $M$ is a faithful finite dimensional
    $A$-module we can find an integer $r>0$ such that $A$ is isomorphic
    to an $A$-submodule of~$\Mr$ as a right $A$-module.  Indeed, if
    $\set{m_1,\dots,m_d}$ is a basis of $M$ then the map
    $A\longrightarrow M^{\oplus d}$ given by $a\mapsto
    (am_1,\dots,am_d)$ is injective. Hence, $A$ is isomorphic to a
    submodule of $\Mr$ whenever $r\ge d=\dim_K M$. Fix such an integer
    $r$. Since $A$ is self-injective the map $A\injection\Mr$ splits,
    so~$\Mr\cong A\oplus N$ as an $A$-module. This completes the proof
    of part~(a).

    Part~(b) follows directly from the definitions since $M$ and $M^{\oplus r}$
    have the same indecomposable summands.
  \end{proof}

  In terms of the representation theory of~$S$, \autoref{P:Equivalence}
  says that there is no harm in assuming that $A$ is a direct summand
  of~$M$ whenever $(A, M)$ is a Schur pair. This will be convenient in
  several of the arguments below.

  The algebra $S$ acts from the left on $M$ as an algebra of
  endomorphisms, with $\phi\cdot m=\phi(m)$, for $\phi\in S$ and $m\in M$.
  Moreover, since $A$ acts faithfully on~$M$, there is an algebra embedding
  $A^\op\injection\End_K(M)$, where we identify $a\in A$ with the
  endomorphism $\rho_a\in\End_K(M)$ given by $\rho_a(m)=ma$, for $m\in M$.
  By construction,
  \[  (\phi\cdot m)a = \phi(m) a = \phi(ma) = \phi\cdot(ma),
            \qquad \text{for all }\phi\in S, m\in M \text{ and } a\in A.
  \]
  Hence, the left and right actions of~$S$ and~$A$ on~$M$
  commute with each other and $A^\op$ is a subalgebra of~$\End_S(M)$.
  The next result will let us show that $A^\op\cong\End_S(M)$.

  \begin{Lemma}\label{L:DoubelCentralizer}
    Suppose that $(A, M)$ is a Schur pair and $r\ge 1$. Then
    $\End_{S_A(M)}(M)\cong\End_{S_A(M^{\oplus r})}(M^{\oplus r})$.
  \end{Lemma}

  \begin{proof} The proof is elementary but we give the argument
    because we need this idea below. There is an embedding
    $\Theta:\End_{S_A(M)}(M)\injection\End_{S_A(M^{\oplus
    r})}(M^{\oplus r})$ given by $\Theta(\phi)=\phi^{\oplus r}$, where
    $\phi^{\oplus r}$ is the endomorphism of $M^{\oplus r}$ given by
    $\phi^{\oplus r}(\underline m)=(\phi(m_1),\dots,\phi(m_r))$, for
    $\underline m=(m_1,\dots,m_r)\in M^{\oplus r}$.  For $1\le i,j\le r$
    let $\pi_i$ be projection onto the $i$th component of~$M^{\oplus r}$
    and let $\sigma_{ij}$ be the map that swaps the $i$th and $j$th
    components.  Then $\pi_i,\sigma_{ij}\in\End_A(M^{\oplus r})$.
    Therefore, if $\phi\in\End_{S_A(M^{\oplus r})}(M^{\oplus r})$ and
    $\underline m\in M^{\oplus r}$ then
    $\phi(\pi_i\underline m)=\pi_i\phi(\underline m)$. Hence, $\phi$
    maps the $i$th component of~$M^{\oplus r}$ into the $i$th component
    so that
    \[
        \phi(\underline m)=\bigl(\phi'_1(\underline m),\dots,\phi'_r(\underline m)\bigr)
                       =\bigl(\phi_1(m_1),\dots\phi_r(m_r)\bigr),
    \]
    for suitable $\phi'_1,\dots,\phi'_r\in\Hom_{S_A(M)}(M^{\oplus r},M)$ and
    $\phi_1,\dots,\phi_r\in\End_{S_A(M)}(M)$. Finally, if $1\le i<j\le r$
    then $\phi(\sigma_{ij}\underline m)=\sigma_{ij}\phi(\underline m)$.
    This implies that $\phi_i(m)=\phi_j(m)$, for all $m\in M$. Hence, the
    map $\Theta$ given above is surjective and the lemma is proved.
  \end{proof}

  The next result shows that $A$ and~$S$ enjoy a \textit{Double
  Centralizer Property} in the sense that $A\cong\End_S(M)^\op$ and
  $S\cong\End_A(M)$. This result is well-known and appears as
  \cite[Theorem~59.6]{C&R}.  We give a self contained proof of this
  result both because it is central to all of the results in this
  section and because the proof is quite short.

  \begin{Theorem}[Double centraliser property]\label{T:DoubleCentraliser}
    Let $(A, M)$ be a Schur pair. Then
    $A^\op\cong\End_S(M)$.
  \end{Theorem}

  \begin{proof}
    By \autoref{P:Equivalence} and \autoref{L:DoubelCentralizer}, we can assume
    that $M\cong A\oplus N$, for some $A$-module~$N$. This reduces the
    proof of the theorem to the claim that
    $\End_S(A\oplus N)\cong A^\op$, where $S=\End_A(A\oplus N)$.
    Recall that $A^\op\cong \End_A({}_AA)$, where ${}_AA$ is the
    left regular representation of~$A$. With this observation the claim now
    follows easily using a similar argument to the proof
    of \autoref{L:DoubelCentralizer}; see \cite[Lemma~59.4]{C&R} for more
    details.
  \end{proof}

  \begin{Corollary}\label{C:SchurFunctor}
    Let $(A, M)$ be a Schur pair such that $A$ is a direct
    summand of~$M$ and let $e\map MA$ be the natural projection map.
    Then $e\in S$, $A\cong eSe$ as an algebra and $M\cong Se$ as an $(S,A)$-bimodule.
  \end{Corollary}

  \begin{proof}
    By definition, $e$ is an endomorphism of $M$ that commutes with the
    action of~$A$ so we can consider~$e$ to be an element of~$S$. Then
    $Se\cong M$ as a left $S$-module. Explicitly, the isomorphism is given by
    $s\mapsto s\cdot 1_A$, where $1_A$ is the identity element of~$A$.
    Hence, by \autoref{T:DoubleCentraliser},
    $A\cong(\End_S(M))^\op\cong (\End_S(Se))^\op\cong eSe$ as an algebra.
    Finally, if we identify $A$ with $eSe$ then $Se\cong M$ as an $(S,A)$-bimodule.
  \end{proof}

  Suppose that $(A, M)$ is a Schur pair. Define a functor $\Fun: \RMod
  S\rightarrow \RMod A$ by $\Fun(X)=X\otimes_{S}M$, for any $X\in S\Mod$
  and $\Fun(g)=g\otimes_{S}\id_{M}$ for any $g\in\Hom_{S}(X,Y)$. Note
  that $\Fun(X)$ is a right $A$-module with $A$-action given by
  $(x\otimes m)a=x\otimes ma$, for all $x\in X, m\in M, a\in A$.
  Moreover, if $A$ is a direct summand of~$M$ with $e\map MA$ be the
  natural projection map, then $\Fun(X)=X\otimes_{S}M\cong
  X\otimes_{S}Se\cong Xe$.

  \begin{Corollary}\label{C:Mprojective}
    The module $M$ is projective as a left $S$-module.
  \end{Corollary}

  \begin{proof}
    In view of \autoref{P:Equivalence}, it is sufficient to consider the
    case when $A$ is a direct summand of~$S$.  Then $M\cong Se$ as an
    $S$-module by \autoref{C:SchurFunctor}, so $M$ is a projective
    $S$-module.
  \end{proof}

  \begin{Corollary} \label{C:FullyFaithful} The functor $\Fun$ is exact
    and fully faithful on projective $S$-modules.
  \end{Corollary}

  \begin{proof}
    The functor $\Fun$ is exact because $M$ is projective as an
    $S$-module by \autoref{C:Mprojective}. To prove that~$\Fun$ is fully
    faithful on projectives it is enough to consider the case when $A$
    is a direct summand of~$M$ by \autoref{P:Equivalence}. Then $A\cong
    eSe$ and $\Fun(X)\cong Xe$ by \autoref{C:SchurFunctor} and the
    discussion above \autoref{C:Mprojective}, for any $S$-module~$X$.
    Hence, $\Fun$ is fully faithful on projectives by
    \autoref{T:DoubleCentraliser}.
  \end{proof}

  By \autoref{C:SchurFunctor}, if $A$ is a direct summand of $M$ then
  $A\cong eSe$ and $\Fun(X)=Xe$, for some idempotent $e\in S$.
  Therefore, in order to understand the connection between $S$-modules
  and $A$-modules we can apply Auslander's theory of quotient functors,
  or Schur functors, as in \cite[\S3.1]{BDK} or \cite{Green}.  Let
  $\set{{L^\lambda}|\lambda\in\P}$ be a complete set of pairwise
  non-isomorphic simple $S$-modules in $\RMod S$, where~$\P$ is an
  indexing set. For each $\lambda\in\P$, define the $A$-module \[
  D^{\lambda}=\Fun(L^{\lambda})\] and set
  $\K=\set{\lambda\in\P|D^\lambda\ne0}$. By the general theory of
  quotient functors, such as in \cite[(6.2g)]{Green},
  $\set{D^\lambda|\lambda\in \K}$ is a complete set of pairwise
  non-isomorphic simple right $A$-modules. By \autoref{P:Equivalence} this is
  true whenever $S=S_A(M)$ and $(A, M)$ is a Schur pair (that is, we do
  not need to assume that $A$ is a direct summand of~$M$).

  \begin{Definition}\label{D:YoungModules}
    Let $\lam\in\P$. Let $P^\lambda$ be the \textbf{projective cover} of $L^\lambda$
    and define the \textbf{Young module} ${Y}^\lambda=\Fun(P^\lambda)$.
  \end{Definition}

  Thus, $P^\lambda$ and $L^\lambda$ are $S$-modules and $D^\lambda$ and
  $Y^\lambda$ are $A$-modules. In the special case when $S$ is the Schur
  algebra and~$A$ is the group algebra of the symmetric group, $Y^\lambda$
  is the usual Young module.

If $E$ is a right $B$-module for an algebra $B$ then the \textbf{socle}
of~$E$, $\soc_B(E)$, is the maximal semisimple (right) submodule of~$E$.
Dually, the \textbf{head} of~$E$, $\hd_B(E)$, is its maximal semisimple (right)
quotient module. When the algebra~$B$ is clear then we simply write $\soc E$
and $\hd E$.

  \begin{Proposition} \label{P:YoungModules} Suppose that $\lam,\mu\in\P$.
    Then the following hold.
    \begin{enumerate}
      \item The Young module $Y^\lambda$ is an indecomposable
      $A$-module.
      \item There is an isomorphism $Y^\lambda\cong Y^\mu$ if and only
      if $\lambda=\mu$.
      \item The Young module $Y^\lambda$ is projective if and only if
      $\lambda\in\K$.
      \item If $\lambda\in\K$ then $Y^\lambda$ is the projective cover
      of~$D^\lambda$.
    \end{enumerate}
  \end{Proposition}

\begin{proof} By \autoref{C:FullyFaithful}, the functor ${\Fun}$
  is fully faithful on projective $S$-modules. Therefore,
  $\End_{A}(Y^\lambda)\cong\End_{S}(P^\lambda)$ is a local $K$-algebra
  since $P^\lambda$ is indecomposable.
  Therefore, $Y^\lambda$ is an indecomposable $A$-module. Further,
  $\Hom_A(Y^\lambda,Y^\mu)\cong\Hom_{S}(P^\lambda,P^\mu)$, so
  $Y^\lambda\cong Y^\mu$ if and only if $\lambda=\mu$. We have
  established parts~(a) and (b), so it remains to prove~(c) and~(d).

  By \autoref{C:SchurFunctor} we may assume that $A$ is a direct summand
  of~$M$, as an $A$-module, so that $A\cong eSe$ for some idempotent
  $e\in S$. Therefore, $\Fun(S)=Se\cong eSe\oplus (1-e)Se$ as right
  $A$-modules.  Therefore, every indecomposable projective right $A$-module is
  isomorphic to $Y^\nu$, for some $\nu\in\P$.

  Fix $\mu\in\K$ and let $Y(\mu)$ be the projective cover of the simple
  $A$-module $D^\mu$.  By the last paragraph, $Y(\mu)\cong Y^\lambda$ for
  some $\lambda\in\P$. Suppose, by way of contradiction, that
  $\lambda\ne\mu$.  As $P^\mu$ is the projective cover of
  $L^\mu$, by applying~$\Fun$ there is a surjective map
  $Y^\mu\surjection D^\mu$ since $\Fun(L^\mu)=D^\mu$. Let
  $s\ge1$ be the multiplicity of $D^\mu$ in the head of~$Y^\mu$.
  As $Y^\mu\not\cong Y^\lambda\cong Y(\mu)$ is indecomposable we can find a
  non-zero map $\theta\map{(Y^\lambda)^{\oplus s}}Y^\mu$ such
  that $D^\mu$ does not appear in the head of~$Y^\mu/\im\theta$.
  Hence, by \autoref{C:FullyFaithful}, there exists a non-zero
  map $\hat\theta\map{(P^\lambda)^{\oplus s}}P^\mu$ such that
  $L^\mu$ does not appear in the head of~$P^\mu/\im\hat\theta$.
  Note that $\im(\hat\theta)\neq P^\mu$ as $\hd(P^\mu)\ncong\hd(P^\lam)$. However, this is a contradiction because $P^\mu$ has simple
  head~$L^\mu$. Hence, $\lambda=\mu$ and $Y(\mu)\cong Y^\mu$ as
  we wanted to show. This completes the proof of parts~(c) and (d).
\end{proof}

\begin{Theorem}\label{T:SDecomposition}
  Suppose that $A$ is a direct summand
  of~$M$ and let $e\map MA$ be the natural projection map.
  \begin{enumerate}
    \item As a right $S$-module,
        $eS\cong\displaystyle\bigoplus_{\mu\in\K}D^\mu\otimes P^\mu$.
    \item As a right $A$-module,
       $M\cong Se\cong\displaystyle\bigoplus_{\lambda\in\P}L^\lambda\otimes Y^\lambda.$
  \end{enumerate}
\end{Theorem}

\begin{proof}
  By \autoref{C:SchurFunctor}, $A\cong eSe$. First consider~(a). There exist non-negative integers $d_\mu$ such that
  \[ eS\cong \bigoplus_{\mu\in\P}(P^\mu)^{\oplus d_\mu}. \]
  Therefore, using \autoref{C:SchurFunctor} for the second equality,
  \[A\cong eSe=\Fun(eS)
     =\bigoplus_{\mu\in\P}\Fun(P^\mu)^{\oplus d_\mu}
     =\bigoplus_{\mu\in\P}(Y^\mu)^{\oplus d_\mu}.
  \]
  The direct summands of~$A$ are necessarily projective
  $A$-modules, so $d_\mu\ne0$ if and only if $\mu\in\K$ and $Y^\mu$ is the projective cover of $D^\mu$ by
  \autoref{P:YoungModules}(d).  Moreover, since $A$ is $K$-split by assumption,
  $d_\mu=\dim D^\mu$ for all $\mu\in\P$. This completes the proof
  of~(a).

  Now consider~(b).
  As a right $S$-module, $S\cong\bigoplus_{\lambda\in\P} L^\lam\otimes P^\lam$
  since $S$ is split over the field~$K$. By definition,
  $Y^\lam=\Fun(P^\lam)$. In view of \autoref{C:SchurFunctor},
  $M\cong Se\cong\Fun(S)\cong\bigoplus_{\lam\in\P} L^\lam\otimes Y^\lam$ as a right
  $A$-module.
\end{proof}

By \autoref{P:YoungModules}, $\set{Y^\mu|\mu\in \K}$ is a
complete set of pairwise non-isomorphic indecomposable projective right
$A$-modules.  As $A$ is self-injective, $\set{Y^\mu|\mu\in \K}$
is also a complete set of pairwise non-isomorphic indecomposable
injective right $A$-modules.

The next few results assume that $A\cong eSe$ as in
\autoref{C:SchurFunctor}. For these results we identify the algebras $A$
and $eSe$ using this isomorphism. For the next result, say that a right
$S$-submodule~$K$ of~$eS$ is \textbf{generated} by $A$
if $K=(K\cap eSe)S$.

\begin{Lemma} \label{L:Ideals}
   Suppose that $A$ is a direct summand of $M$ and let $e\map MA$ be the natural projection map. Then:
  \begin{enumerate}
   \item The map $I\mapsto IS$ defines an inclusion preserving
    bijection between the set of right $A$-submodules of $eSe$ and
    the set of right $S$-submodules of~$eS$ generated by~$A$.
    \item If $I$ and $J$ are right $A$-submodules of $eSe$ then
        $\Hom_A(I,J)\cong\Hom_S(IS,JS).$
    \item If $I$ and $J$ are right $A$-submodules of $eSe$ then $I\cong J$
    as $A$-modules if and only if $IS\cong JS$ as $S$-modules.
  \end{enumerate}
\end{Lemma}

\begin{proof} Let $I$ be a right $A$-submodule of $eSe$. We claim that
  $I=IS\cap eSe$. Certainly,
  $I\subseteq IS\cap eSe$. Conversely, if $x\in IS\cap eSe$ then $x=xe$ and
  we can write $x=\sum_{s\in S}a_s s$, for some $a_s\in I$. Therefore,
  \[ x=xe = \sum_{s\in S}a_s se=\sum_{s\in S}a_s (ese)\in I, \]
  since $ese\in eSe$ and $I$ is a right $A$-module. Therefore,
  $I=IS\cap eSe$ as claimed. Hence, the ideal $IS$ is generated by~$A$.
  Moreover, if $I$ and $J$ are right $A$-submodules of $eSe$ then $IS=JS$ if
  and only if $I=J$. We have now proved part~(a).

  As (c) follows immediately from~(b), it remains to prove~(b). Let $I$ and $J$ be $A$-submodules of $eSe$.
  By the last paragraph, any homomorphism $\psi\map{IS}{JS}$ restricts
  to a well-defined $A$-module homomorphism $\psi_e\map IJ$ since
  $A\cong eSe$. Conversely, since $A$ is self-injective any
  homomorphism between two right ideals of $A$ is given by left
  multiplication by \cite[Theorem~61.2]{C&R}. That is, if $\psi_e\map
  IJ$ is an $A$-module homomorphism then there exists an $a\in eSe$ such
  that $\psi(x)=ax$, for all $x\in I$. Therefore, there is a
  well-defined $S$-module homomorphism $\psi\map{IS}{JS}$ given by
  $\psi(y)=ay$, for all $y\in IS$. By construction, $\psi$ is uniquely
  determined by~$\psi_e$, and it restricts to~$\psi_e$, so~(b) follows.
\end{proof}

Suppose that the algebra $A$ comes equipped with an
anti-involution~$\ast$. If $X$ is any right $A$-module then the dual
module $X^\ast=\Hom_K(X,K)$ becomes a right $A$-module with action
$(fa)(x) = f(xa^\ast)$, for $a\in A$, $f\in X^\ast$ and $x\in X$.

\begin{Definition}\label{D:SelfDualSchurPair}
  A \textbf{self-dual Schur pair} is a Schur pair $(A, M)$ where $A$ is
  equipped with an anti-involution $\ast$ such that $M\cong M^*$ as $A$-modules.
\end{Definition}

Except for the last two results in this section, we now assume that $(A, M)$
is a self-dual Schur pair and we fix an $A$-module isomorphism $\theta:
M\bijection M^*$. Following \cite[(1.5)]{DJ:qWeyl}, define an
anti-isomorphism
\[ \tau\map{\End_A(M)}\End_A(M^*);\quad s\mapsto s^\tau , \]
where $s^\tau(f)=f\circ s$, for $s\in S=\End_A(M)$ and $f\in M^*$.
Then the anti-isomorphism~$\ast$ on~$A$ induces an anti-isomorphism
$\widetilde{\tau}$ on~$S=\End_A(M)$ that is given by
$\widetilde{\tau}(s):=s^*:=\theta^{-1}\circ s^{\tau}\circ\theta$, for
all $s\in S$. By definition, if $s\in S$ and $m\in M$ then
$s^*m=\theta^{-1}(\theta(m)\circ s)$ so that
$\theta(s^*m)=\theta(m)\circ s$ or, equivalently,
$\theta(sm)=\theta(m)\circ s^*$. Hence, $\theta$ is also an $S$-module
isomorphism.

\begin{Lemma}\label{L:MSselfdual}
 Suppose that $(A, M)$ is a self-dual Schur pair. Then $M$ is self-dual
 and projective-injective as an $S$-module.
\end{Lemma}

\begin{proof}
  In the last paragraph we observed that $\theta$ is an $S$-module
  homomorphism, so $M\cong M^*$ as an $(S,A)$-bimodule. In particular,
  $M$ is self-dual as an $S$-module. By \autoref{C:Mprojective}, $M$ is
  projective as an $S$-module so $M^*\cong M$ is injective giving the
  remaining claim.
\end{proof}

\begin{Definition} \label{D:rightRdfn} Suppose that $(A, M)$ is a self-dual Schur pair such that $A$ is a direct
  summand of $M$ and let $e\map MA$ be the natural projection map. Assume that $e^\ast=e$. For any $(S,A)$-bimodule $X$
let $X_R$ be the $(A,S)$-module such that $X_R=X$ as a $K$-vector space and
$axs=s^{\ast}xa^\ast$, for $a\in A, s\in S$ and $x\in M$.
\end{Definition}

Suppose that $A$ is a direct summand of $M$ and let $e\map MA$ be the
natural projection map. By definition and \autoref{C:SchurFunctor},
$S_A(M)=\End_A(M)\cong\End_A\bigl(Se\bigr)$. The next lemma, which
is part of the motivation for introducing self-dual Schur pairs, allows
us to replace~$Se$ with~$eS$. Note that $eS$ is a left $A$-module since
$A\cong eSe$.

\begin{Lemma}\label{L:ee*}
  Suppose that $(A, M)$ is a self-dual Schur pair such that $A$ is a
  direct summand of $M$ and let $e\map MA$ be the natural projection
  map. Assume that $e^\ast=e$. Then $\ast$ induces an isomorphism
\[
    (S_A(M))^{\op}\cong S_A(M),
\]
and right multiplication induces an algebra isomorphism $S_A(M)^{\op}\cong\End_{A}(eS_A(M))$.
\end{Lemma}

\begin{proof} We identify $A$ with $eS_A(M)e$. The first isomorphism
  follows from the induced anti-isomorphism $\ast$ of $S_A(M)$.
 Since $e=e^*$, and for $a\in A, s\in S$, the map $es\mapsto s^\ast e$
 defines an $(A,S)$-bimodule isomorphism $eS\cong (Se)_R=M_R$ (see
 \autoref{D:rightRdfn} for the definition of $M_R$). By definition, $M$
 is a faithful left $S_A(M)$-module. It follows that $eS_A(M)$ is a
 faithful right $S_A(M)$-module. Hence, right multiplication induces
 an injective algebra homomorphism
  $S_A(M)^{\op}\hookrightarrow\End_{A}(eS_A(M))$. On the other hand,
  \[
  \dim \End_{A}(eS_A(M))=\dim \End_{A}(S_A(M)e)=\dim S_A(M).
  \]
  It follows that the injection above is an isomorphism.
\end{proof}

\begin{Lemma} \label{L:KeyStep}
   Suppose that $(A, M)$ is a self-dual Schur pair such that $A$ is a direct
   summand of $M$ and let $e\map MA$ be the natural projection map. Assume that $e^\ast=e$. Let
   $X$ be a non-zero right $S$-submodule of $eS$. Then $X\cap eSe$ is a
   non-zero right $A$-submodule of~$eSe$.
\end{Lemma}

\begin{proof} It is easy to see that $X\cap eSe$ is a right $A$-submodule
  of $eSe$. We need to prove that $X\cap eSe\ne\{0\}$.

  Before we start the proof observe that we can regard $eS$ as a
  left $A$-module because $A\cong eSe$. Let $D=\soc_A(eS)$ and let
  $Q=Q_D$ be the injective hull of~$D$. As $A$ is self-injective
  $Q$ is also projective as a left $A$-module. Therefore, we can find
  an integer $t\ge0$ such that $Q$ embeds into $A^{\oplus t}$ as
  a left $A$-module. Recalling that $Q$ is the injective hull of
  $D=\soc_A(eS)$, this implies that there exists a left $A$-module
  homomorphism $\theta\map{eS}A^{\oplus t}$ such that $\theta$
  restricts to the identity map on~$D$.

  Now are now ready to prove the lemma. Recall that $\widetilde{\tau}(s)=s^\ast$ for each $s\in S$. By assumption $\widetilde{\tau}(e)=e^\ast=e$. Since $A_{A}$ is a direct summand of $Se$ as a right $A$-module. It follows that the left regular $A$-module ($\cong\widetilde{\tau}(A_A)$) is a direct summand of $eS=\widetilde{\tau}(Se)$ as a left $A$-module. Fix a non-zero element $x\in X$.
  Then $x=ex\in eS$ so we can find a non-zero element $a=eae\in eSe\cong
  A$ such that $ax\in D=\soc_A(eS)$ because $0\neq
  \soc_A(eSex)\subseteq\soc_{A}(eS)=D$. Therefore, $\theta(ax)\ne0$. By
  composing~$\theta$ with a suitable projection of~$A^{\oplus t}$ onto~$A$, and then using the inclusion $A\injection eS$, we obtain
  a left $A$-module homomorphism $\vartheta\map{eS}{eS}$ such that
  $\im\vartheta\subseteq A$ and $\vartheta(ax)\ne0$. Consequently,
  $\vartheta(x)\ne0$ since $\vartheta(ax)=a\vartheta(x)$. Write
  $\vartheta(x)=ye+z(1-e)$, for some $y,z\in eS$. Then
  $aye+az(1-e)=\vartheta(ax)\in eSe$, so $aye\ne0=az(1-e)$. In
  particular, $ye=eye$ is a non-zero element of~$A$.
  To complete the proof recall that $S^{\op}\cong\End_A(eS)$ by
  \autoref{L:ee*}. Therefore, we may assume that
  $\vartheta\in S$ and consider the endomorphism $e\circ\vartheta$
  of~$Se$. Since $X$ is an $S$-module,
  $ye=\vartheta(x)e=x(\vartheta^\ast\circ e)\in X\cap eSe$, where $x\vartheta^\ast=\vartheta(x)$ follows by definition. That is, $ye$ is a non-zero
  element of $X\cap eSe$, so $X\cap eSe\ne0$ as we wanted to show.
\end{proof}

We can now prove the key result of this section.  Recall that
$\set{L^\lambda|\lambda\in\P}$ is a complete set of irreducible
$S$-modules and that $\set{D^\mu|\mu\in\K}$ is a complete set of
irreducible $A$-modules.

\begin{Theorem} \label{T:Socle}
   Suppose that $(A,M)$ is a self-dual Schur pair such that $A$ is a direct summand
   of $M$ and let $e\map MA$ be the natural projection map. Assume that $e^\ast=e$. Let
   $\mu\in\P$. The simple $S$-module $L^\mu$ is isomorphic to a
   submodule of $eS$ if and only if~$\mu\in\K$.
\end{Theorem}

\begin{proof} Suppose that $L$ is an irreducible right $S$-submodule of~$eS$. Then
  $D=L\cap eSe\ne0$ by \autoref{L:KeyStep}. Therefore, $DS$ is a
  non-zero $S$-submodule of~$L$ so that $L=DS$ since $L$ is
  irreducible. In view of \autoref{L:Ideals}(a), $D$ is an
  irreducible right $A$-submodule of~$eSe$. Moreover, since
  $xe=x$ for all $x\in D$ it follows that $Le\ne0$ since
  $D\subseteq Le$. Therefore, $Le=D$ since $D$ is a non-zero
  $A$-submodule of~$Le$ and $Le$ is irreducible. That is,
  $\Fun(L)=Le=D$. Consequently, if $L\cong L^\mu$ then
  $D\cong\Fun(L^\mu)=D^\mu$ so that $\mu\in\K$.

  Conversely, because $A$ is self-injective, if $\mu\in\K$ then  there
  exists an irreducible submodule~$D$ of~$eSe$ such that $D\cong D^\mu$
  as a right $A$-module. By \autoref{L:Ideals}(a), $L=DS$ is an irreducible
  right $S$-submodule of~$eS$. By the argument of the last
  paragraph $\Fun(L)\cong D^\mu$. Hence, $L\cong L^\mu$ and $eS$
  has a submodule isomorphic to~$L^\mu$ for all~$\mu\in\K$.
\end{proof}

Note that \autoref{T:Socle} would follow from \autoref{T:SDecomposition}
if $eS$ were self-dual as an $S$-module.

Recall that $\soc_S(eS)$ is the socle of~$eS$ as an $S$-module.
\autoref{T:Socle} says that $L^\mu$ appears in $\soc_S(eS)$ if
and only if $\mu\in\K$.

\begin{Corollary} \label{C:Socle} Let $\lam,\mu\in\P$. Then $L^\mu$
  appears in the socle of $P^\lambda$ only if $\mu\in \K$.
\end{Corollary}

\begin{proof} By definition, $eS$ is a faithful right $S$-module.
  Therefore, as in the proof of \autoref{P:Equivalence}, $S$ is isomorphic to a
  submodule of $(eS)^{\oplus t}$, for some $t\ge1$. It follows that $P^\lambda$ is
  isomorphic to a submodule of~$(eS)^{\oplus t}$. In
  particular, the socle of $P^\lambda$ is isomorphic to a submodule of
  $\soc((eS)^{\oplus t})=(\soc eS)^{\oplus t}$. Hence,~$L^\mu$
  appears in $\soc(P^\lambda)$ only if $\mu\in\K$ by
  \autoref{T:Socle}.
\end{proof}

A \textbf{block} of an algebra~$B$ is an indecomposable two-sided ideal. Every algebra
  has a unique decomposition $B=B_1\oplus\dots\oplus B_z$ into a direct
  sum of blocks. A $B$-module $X$
  \textbf{belongs} to the block $B_r$ if $X=XB_r$. It is easy to see
  that indecomposable modules belong to a unique block. In particular,
  the block decomposition of~$B$ induces an equivalence relation on the
  simple $B$-modules where two simples are equivalent if they belong to
  the same block. As is well-known, this equivalence relation  coincides
  with the \textbf{linkage classes} of simple $B$-modules, where the
  equivalence relation for the linkage classes is generated by $D\sim E$
  if $\Hom_B(P_D,P_E)\ne0$, where $P_D$ and $P_E$ are the projective
  covers of the simple modules $D$ and~$E$.

The following corollary is well-known. We omit the proof following the suggestion of the referee.

  \begin{Corollary}\label{C:Blocks} Suppose that $(A,M)$ is a self-dual Schur pair.
    The functor $\Fun$ induces a one-to-one correspondence between the
    blocks of~$S$ and the blocks of~$A$.
  \end{Corollary}

\begin{Lemma}\label{L:YSelfDual}
  Suppose that $A$ is a symmetric algebra and $(D^\mu)^*\cong D^\mu$, for all $\mu\in\K$. Then $(Y^\mu)^*\cong Y^\mu$ as $A$-modules for all $\mu\in\K$.
\end{Lemma}

\begin{proof} By  \cite[Theorem~3.1 (d)]{Rickard:SymmetricAlgs}, $\Hom_A(Y^\lam,D^\lam)\cong\bigl(\Hom_A(D^\lam, Y^\lam)\bigr)^\ast$ for  all $\lam\in\K$.
Since $(D^\mu)^*\cong D^\mu$ for all $\mu\in\K$, the lemma follows at once. Alternatively, a direct
  proof can be given by adapting the argument of \cite[Theorem~6]{Alperin}.
\end{proof}

From now until \autoref{C:KleshSocle}, we assume that  $(A, M)$ is a
self-dual Schur pair such that $A$ is a direct summand of $M$ and that
$A$ is a symmetric algebra such that $(D^\mu)^*\cong D^\mu$, for all
$\mu\in\K$.  Let $e\map MA$ be the natural projection map and assume
that $e^\ast=e$.

If $N\in\RMod{S}$ then $\Fun(N)=Ne$ by definition. It is easy to check
that the map $fe\mapsto\psi(fe): xe\mapsto f(xe)$, $\forall\, f\in
N^\ast, x\in N$ defines a right $A$-module isomorphism $\psi: N^\ast
e\cong (Ne)^\ast$. In other words, the Schur functor $\Fun$ commutes
with the duality functor . As a consequence, it follows easily
that
\begin{equation}\label{selfdual2}
   (L^\mu)^*\cong L^\mu, \qquad\text{for all }\mu\in\K.
\end{equation}

The next result shows that the projective-injective $S$-modules are
indexed by the simple $A$-modules. The assumption that the simple
$A$-modules are self-dual is used to ensure that the head and socle
of~$P^\lambda$ are isomorphic whenever $P^\lambda$ is self-dual. In
applications $A$ is usually a cellular algebra, in the sense
of~\cite{GL}, in which case this assumption is automatic.

\begin{Theorem}\label{T:Equivalence}
  Suppose that $(A, M)$ is a self-dual Schur pair such that $A$ is a
  direct summand of $M$ and let $e\map MA$ be the natural projection
  map. Assume that $e^\ast=e$ and $(D^\mu)^*\cong D^\mu$, for  all
  $\mu\in\K$.  Let $\lambda\in\P$. Then the following are equivalent:
  \begin{enumerate}
    \item $\lambda\in\K$,
    \item $D^\lambda\ne0$,
    \item $L^\lambda$ is a right $S$-submodule of $M_R$,
    \item $P^\lambda$ is a direct summand of $M_R$ as a right $S$-module,
    \item $P^\lambda$ is a projective-injective right $S$-module,
    \item $P^\lambda$ is self-dual.
  \end{enumerate}
\end{Theorem}

\begin{proof}
  By construction, $M_R\cong eS$ and, by definition, $\lambda\in\K$ if
  and only if $D^\lambda\ne0$. Further, (a) and (c) are equivalent by
  \autoref{T:Socle} and (a) and (d) are equivalent by
  \autoref{T:SDecomposition}. Hence, (a), (b), (c) and (d) are
  equivalent. To complete the proof we show that (a) $\implies$ (f)
  $\implies$ (e) $\implies$ (d).

  First suppose that (a) holds so that $\lambda\in\K$. By definition,
  $P^\lambda$ has a simple head $L^\lambda$. Then $P^\lambda$ is an
  indecomposable direct summand of~$M_R$ by \autoref{T:SDecomposition}.
  By \autoref{L:MSselfdual}, $M^\ast\cong M$ as a left $S$-module. It follows from \autoref{D:rightRdfn} that
  $(M_R)^\ast\cong M_R$ as a right $S$-module. Therefore, $(P^\lambda)^*$ is an indecomposable summand of~$M_R$, so $(P^\lambda)^*\cong P^\mu$ for some
  $\mu\in\K$. Hence, by (\ref{selfdual2}), $P^\mu$ has
  simple socle~$L^\lam$ which implies that there exists a non-zero
  homomorphism~$f$ from $P^{\lam}$ to $P^\mu$ such that $\im(f)\cong
  L^\lambda$.  As the functor $\Fun$ is fully faithful there exists a
  non-zero homomorphism $f_Y$ from $Y^{\lam}$ to~$Y^\mu$ such that
  $\im(f_Y)\cong D^\lambda\neq 0$. By \autoref{L:YSelfDual}, $Y^\mu$ is
  self-dual so it has simple socle $D^\mu$.  Therefore, $\lam=\mu$ so
  that $(P^\lambda)^*\cong P^\lambda$ is self-dual. Hence, (a)
  $\implies$ (f).

  Now suppose that (f) holds so that $P^\lambda$ is self-dual. By
  definition, $P^\lambda\cong(P^\lambda)^*$ is projective so it is also
  injective and (e) holds.

  Finally, suppose that (e) holds so that $P^\lambda$ is a
  projective-injective $S$-module. Since $M_R$ is a faithful right $S$-module, $S$
  embeds into $M_R^{\oplus t}$, for some $t\ge0$, by the argument of
  \autoref{P:Equivalence}. Fix such a~$t$. Since
  $(P^\lambda)^*$ is projective it embeds into $M_R^{\oplus t}$. Taking
  duals there is a surjection $M_R^{\oplus t}\cong(M_R^*)^{\oplus
  t}\surjection P^\lambda$. As $P^\lambda$ is projective this map
  splits, so $P^\lambda$ is isomorphic to an indecomposable submodule
  of~$M_R^{\oplus t}$. Hence, $P^\lambda$ is isomorphic to an
  indecomposable submodule of~$M_R$.  Hence, (e) $\implies$ (d) as we
  wanted to show.
\end{proof}

If $S$ is a quasi-hereditary algebra then \autoref{T:Equivalence}(e) says
that for any $\mu\in\K$, the projective indecomposable module $P^\mu$ is also the injective
hull of~$L^\mu$ and the indecomposable (partial) tilting module
corresponding to~$\mu$; see \autoref{S:SchurAlgebras}.

\begin{Corollary} \label{C:KleshSocle} Suppose that $X$ is
  a right $S$-module such that $[\soc X:L^\blam]\ne0$ only if
  $\lambda\in\K$. Then \[[\soc X:L^\lambda]=[\soc
  \Fun(X):D^\lambda]\qquad \text{ for all }\blam\in\K.\] In particular,
  if $P$ is a projective $S$-module then $[\soc P:L^\lambda]=[\soc
  \Fun(P):D^\lambda]$, for all $\blam\in\K$.
\end{Corollary}

\begin{proof}If $X=P$ is projective then $[\soc P:L^\lambda]\ne0$ only
  if $\lambda\in\K$ by \autoref{C:Socle}. Suppose, then, that $X$ is a right
  $S$-module such that  $[\soc X:L^\blam]\ne0$ only if $\lambda\in\K$.
  By exactness, $\Fun(\soc X)\subseteq\soc\Fun(X)$. We need to show that
  the reverse inclusion $\soc\Fun(X)\subseteq \Fun(\soc X)$ holds. By
  assumption,
  \[ \soc X=\bigoplus_{\bmu\in \K}(L^\mu)^{\oplus x_\mu} \]
  for some non-negative integers $x_\mu$. Let $I_X$ be the injective hull
  of $\soc X$. Then $I_X=\bigoplus_{\mu\in\K}(P^\mu)^{\oplus x_\mu}$
  since~$P^\mu$ is self-dual for all $\mu\in\K$ by
  \autoref{T:Equivalence}. Therefore we have injections
  \[\bigoplus_{\mu\in\K}(L^\mu)^{\oplus x_\mu}=\soc X
       \injection X \injection
       I_X=\bigoplus_{\mu\in\K}(P^\mu)^{\oplus x_\mu}.
  \]
  Applying the functor $\Fun$ gives injections
  \[\bigoplus_{\mu\in\K}(D^\mu)^{\oplus x_\mu}\injection
       \Fun(\soc X) \injection\Fun(X)\injection
       \Fun(I_X)=\bigoplus_{\mu\in\K}(Y^\mu)^{\oplus x_\mu}.
  \]
  In particular, the socle of $\Fun(X)$ is contained in
  $\soc\Fun(I_X)=\bigoplus_\mu(D^\mu)^{\oplus x_\mu}=\Fun(\soc X)$.
  Hence, $\soc\Fun(X)=\Fun(\soc X)$  and the corollary follows.
\end{proof}

We conclude this section by relating hom-spaces of certain $A$-modules
and $S$-modules. These results do not assume that $(A,M)$ is a self-dual
Schur pair. If $B$ is an algebra and $X$ is a subset of a left or
right $B$-module, respectively, define the left and right
\textbf{annihilators} of~$X$ to be
\begin{align*}
  \L_B(X) & =\set{b\in B|bx=0\text{ for all }x\in X}\\
  \R_B(X) & =\set{b\in B|xb=0\text{ for all }x\in X}.
\end{align*}
Set $\LR_B(X)=\L_B(\R_B(X))$ and $\RL_B(X)=\R_B(\L_B(X))$. If
$X=\{x\}$ write $\L_B(x)=\L_B(X)$ and similarly for
$\R_B(x)$, $\LR_B(x)$ and~$\RL_B(x)$. In particular,
$\LR_B(x) = \set{b\in B|ba=0\text{ whenever }xa=0, \text{ for }a\in B}$.

The key property of the double annihilators $\LR_B(x)$ and $\RL_B(x)$
that we need is the following well-known fact, which  characterises
self-injective algebras.

\begin{Lemma}[{\cite[Theorem~61.2]{C&R}}]\label{L:DoubleAnnihilators}
  Suppose that $B$ is a self-injective algebra and that $x\in B$. Then
  $\LR_B(x)=Bx$ and $\RL_B(x)=xB$.
\end{Lemma}


\begin{Lemma}\label{L:Annihilators} Let $M$ be an $(S,A)$-bimodule and $m\in M$.
  \begin{enumerate}
    \item If $X$ is a left $S$-submodule of $S$
    then $\Hom_{S}(S m,X)\cong \RL_S(m)\cap X$ as vector spaces.
    \item If $Y$ is a right $A$-submodule of $A$
    then $\Hom_{A}(mA,Y)\cong \LR_A(m)\cap Y$ as vector spaces.
  \end{enumerate}
\end{Lemma}

\begin{proof}
  (a) If $x\in\RL_S(m)\cap X$ then the map $sm\mapsto sx$ is an $S$-module
  homomorphism. Conversely, if $f\in\Hom_S(Sm,X)$ then $f(m)\in\RL_S(m)\cap X$
  since $sf(m)=f(sm)=0$ whenever $s\in\L_S(m)$.  It is straightforward
  to check that these maps are mutually inverse isomorphisms.

  (b) By \autoref{L:DoubleAnnihilators}, it is enough to show that
  $\Hom_A(mA, Y)\cong\LR_A(m)\cap Y$. The argument to prove this is
  almost identical to part~(a) in that if $x\in\LR_A(m)\cap Y$ then the
  map $ma\mapsto xa$, for $a\in A$, is an $A$-module homomorphism.
  Conversely, let $f\in\Hom_{A}(mA,Y)$ and $x=f(m)$. We need to show
  that $xa=0$ whenever $a\in\R_A(m)$. That is, $xa=0$ whenever $ma=0$.
  In fact, $xa=f(m)a=f(ma)=0$.
\end{proof}

If $A$ is self-injective and $m\in A$ in \autoref{L:Annihilators}(b) then
$\LR_A(m)=Am$ by \autoref{L:DoubleAnnihilators}, so  $\Hom_{A}(mA,Y)\cong Am\cap Y$.

\section{Cyclotomic Hecke algebras}\label{S:HeckeAlgebras}

The main results of this paper apply the results of \autoref{S:SchurPairs}
to the cyclotomic Schur algebras~$\Sn$, which were introduced
in~\cite{DJM:cyc}. In order to do this we first need to produce a Schur
pair $(A, M)$ that defines the cyclotomic Schur algebras. Using the notation of
\autoref{D:SchurPair}, the algebra $A$ will be a cyclotomic Hecke algebra
$\Hn$ and $M$ will be a direct sum of permutation modules for~$\Hn$.
This section recalls the representation theory of
the cyclotomic Hecke algebras of type~$A$ and proves some new results
about the simple modules of these algebras that will be needed later.

Fix a non-negative integer $n$. Let $\Sym_n$ be the symmetric group on~$n$
letters and let $\set{s_1,\dots,s_{n-1}}$ be the standard set of Coxeter
generators for~$\Sym_n$, where $s_r=(r,r+1)$ for $1\le r<n$. If
$w\in\Sym_n$ then the \textbf{length} of $w$ is
$\ell(w)=\min\set{k\ge0|w=s_{r_1}\dots s_{r_k}}$.

Fix a field $K$ and $\xi\in K^\times$. If $k\in\Z$ then the
\textbf{$\xi$-quantum integer} is
\begin{equation}\label{E:QInt}
[k]=[k]_\xi = \begin{cases}
  \phantom{-}\xi+\xi^3+\dots+\xi^{2k-1},&\text{if }k\ge0,\\
            -\xi^{-1}-\xi^{-3}-\dots-\xi^{2k+1},&\text{if }k<0.
\end{cases}
\end{equation}
So, $[k]_\xi=-[-k]_{\xi^{-1}}$ and if $\xi^2\ne1$ then
$[k]=(\xi^{2k}-1)/(\xi-\xi^{-1})$.

The \textbf{quantum characteristic} of $\xi$ is the smallest positive
integer~$e$ such that $[e]_\xi=0$, where we set~$e=\infty$ if $[k]\ne0$ for
all $k\in\N$. Notice that $\xi$ and $\xi^{-1}$ have the same quantum
characteristic and if~$\xi=1$ then the quantum characteristic
of~$\xi$ is the characteristic of~$K$.

Finally fix an integer $\ell>0$ and a \textbf{multicharge}
$\charge=(\kappa_1,\dots,\kappa_\ell)\in\Z^\ell$.

\begin{Definition}[Cyclotomic Hecke algebras of
  type~$A$\cite{AK,HuMathas:SeminormalQuiver}]\label{D:Hecke}
  The \textbf{cyclotomic Hecke algebra} of type~$A$ with Hecke
  parameter $\xi$ and multicharge~$\charge$ is the unital associative
  $K$-algebra $\Hn=\Hn^\charge(\xi)$ with generators
  $T_1,\dots,T_{n-1}$, $L_1,\dots,L_n$ and relations
  \begin{align*}
      \textstyle\prod_{l=1}^\ell(L_1-[\kappa_l])&=0,  &
      (T_r-\xi)(T_r+\xi^{-1})&=0,  \\
      L_rL_t&=L_tL_r, &
    T_rT_s&=T_sT_r &\text{if }|r-s|>1,\\
    T_sT_{s+1}T_s&=T_{s+1}T_sT_{s+1}, &
    T_rL_t&=L_tT_r,&\text{if }t\ne r,r+1,\\
    \span\span L_{r+1}=T_rL_rT_r+T_r,\span\span\span
  \end{align*}
where $1\le r<n$, $1\le s<n-1$ and $1\le t\le n$.
\end{Definition}

These algebras are almost the same as the Ariki-Koike algebras
introduced in~\cite{AK} except that the presentation of
\autoref{D:Hecke} changes the algebras when $\xi^2=1$. This allows the
so-called \textit{degenerate} ($\xi^2=1$) and \textit{non-degenerate} ($\xi^2\ne1$)
cases, to be treated simultaneously.
See~\cite[\S2]{HuMathas:SeminormalQuiver} for more details.

\begin{Remark}
   \autoref{D:Hecke} is a renormalisation of the presentation of the
   cyclotomic Hecke algebras given in~\cite{HuMathas:SeminormalQuiver}.
   Explicitly, if $\tilde T_r$ and $\tilde L_s$ are the generators of
   the algebra $\tilde\HH_n(\xi^2,\charge)$ given in
   \cite[Definition~2.2]{HuMathas:SeminormalQuiver} then
   $T_r=\xi^{-1}\tilde T_r$ and $L_s=\xi^{-1}\tilde L_s$, for
   $1\le r<n$ and $1\le s\le n$. (In the notation
   of~\cite{HuMathas:SeminormalQuiver}, the cyclotomic parameters of
   $\tilde\HH_n(\xi^2,\charge)$ are $Q_r=\xi^{2\kappa_r}$, for
   $1\le r\le\ell$.)
\end{Remark}

If $w\in\Sym_n$ has reduced expression $w=s_{r_1}\dots s_{r_k}$, so that
$k=\ell(w)$, then set $T_w=T_{r_1}\dots T_{r_k}$.  Then~$T_w$ depends
only on~$w$ and
$
   \set{L_1^{a_1}\dots L_n^{a_n}T_w|0\le a_r<\ell\text{ and } w\in\Sym_n}
$
is a $K$-basis of~$\Hn$ by the argument of~\cite{AK}. Let $*$ be the
unique anti-isomorphism of~$\Hn$ that fixes each of the generators.
Then~$T_w^*=T_{w^{-1}}$ and $L_m^*=L_m$, for all $w\in\Sym_n$ and
$1\le m\le n$.

\begin{Theorem}[%
  Malle-Mathas~\cite{MM:trace} and Brundan~\cite{Brundan:degenCentre}%
  ]\label{T:SymmetricAlgebra}
  The algebra $\Hn$ is a symmetric algebra with non-degenerate trace
  form~$\tau$. In particular, $\Hn$ is a self-injective algebra.
\end{Theorem}

The paper~\cite{MM:trace} proves this result when $\xi^2\ne1$ whereas
\cite{Brundan:degenCentre} treats the case when $\xi^2=1$. In both cases
the trace form~$\tau$ is described explicitly, however, dual bases are
known only when $\ell=1,2$.

A \textbf{partition} of an integer $m$ is a weakly decreasing sequence
$\lambda=(\lambda_1\ge\lambda_2\ge\dots)$ of non-negative integers such
that $|\lambda|=\sum_i\lambda_i=m$. A \textbf{multipartition}, or
$\ell$-partition, of~$n$ is an ordered $\ell$-tuple
$\blam=(\lambda^{(1)}|\dots|\lambda^{(\ell)})$ of partitions such that
$|\blam|=|\lambda^{(1)}|+\dots+|\lambda^{(\ell)}|=n$.  The
\textbf{diagram} of $\blam$ is the set
\[[\blam]=\set{(k,r,c)|r\ge 1, 1\le c\le\lambda^{(k)}_r \text{ and }1\le k\le \ell}.\]
A multipartition is uniquely determined by its diagram. A \textbf{node}
is any element of the diagram of some multipartition.

Let $\Parts$ be
the set of multipartitions of $n$. If $\blam,\bmu\in\Parts$
then $\blam$ \textbf{dominates} $\bmu$, written~$\blam\gedom\bmu$,
if
\[
   \sum_{t=1}^{s-1}|\lambda^{(t)}|+\sum_{i=1}^k\lambda^{(s)}_i
          \ge\sum_{t=1}^{s-1}|\mu^{(t)}|+\sum_{i=1}^k\mu^{(s)}_i,
          \qquad
          \text{for all }1\le s\le \ell \text{ and } k\ge1.
\]
Dominance defines a partial order on $\Parts$.

If $X$ is a set then an $X$-valued $\blam$-tableau is a function
$\bT\map{[\blam]}X$. If $\bT$ is a $\blam$-tableau write
$\Shape(\bT)=\blam$. For convenience we identify
$\bT=(\bT^{(1)},\dots,\bT^{(\ell)})$ with a labeling of the diagram
$[\blam]$ by elements of $X$ in the obvious way. Thus, we can talk of
the components, rows and columns of $\bT$.

A \textbf{standard $\blam$-tableau} is a map
$\t\map{[\blam]}\{1,2,\dots,n\}$ such that for $s=1,\dots,\ell$ the
entries in each row of $\t^{(s)}$ increase from left to right and the
entries in each column of $\t^{(s)}$ increase from top to bottom.  Let
$\Std(\blam)$ be the set of standard $\blam$-tableaux and set
\[ \Std(\Parts)=\bigcup_{\blam\in\Parts}\Std(\blam)
\quad\text{and}\quad
\Std^2(\Parts)=\bigcup_{\blam\in\Parts}\Std(\blam)\times\Std(\blam).
\]
If $\t\in\Std(\Parts)$ and $1\le m\le n$ let $\comp_m(\t)=k$ if~$m$
appears in the $k$-th component of~$\t$.

The \textbf{conjugate} of~$\blam$ is the multipartition~$\blam'$ whose
diagram is obtained from $[\blam]$ by reversing the order of the components
and then swapping rows and columns.  Thus,
$[\blam']=\set{(k,r,c)|(\ell-k+1,c,r)\in[\blam]}$.  The
\textbf{conjugate} of $\t$ is the $\blam'$-tableau~$\t'$ obtained by
reversing the order of the components in~$\t$ and then transposing the
tableau in each component (that is, swapping rows and columns).

Let $\tlam$ be the standard $\blam$-tableau such that the numbers
$1,2,\cdots,n$ are entered in order from left to right along the rows of
$\t^{\lambda^{(1)}}$, and then $\t^{\lambda^{(2)}}, \dots,
\t^{\lambda^{(\ell)}}$. Similarly, let $\tllam$ be the standard
$\blam$-tableau with the numbers $1,2,\cdots,n$ entered in order down
the columns of $\t_{\blam^{(\ell)}},\dots,\t_{\blam^{(1)}}$. Then
$(\tlam)'=\t_{\blam'}$, for all $\blam\in\Parts$.

If $\t$ is a standard $\blam$-tableau let $d(\t),d'(\t)\in\Sym_n$ be the
unique permutations such that $\t=\tlam d(\t)$ and $\t=\tllam d'(\t)$.
Let $w_\blam=d(\tllam)$. It is easy to see that $d(\t')=d'(\t)$
and $w_\blam=d(\t)d(\t')^{-1}$, for all $\t\in\Std(\blam)$.

If $\t$ is a tableau with entries in $\set{1,2\dots,n}$ let $\Row(\t)$
and $\Col(\t)$ be the subgroups of~$\Sym_n$ that stabilise the rows and
columns of~$\t$, respectively.  If $\blam\in\Parts$ let
$\Sym_\blam=\Sym_{\lambda^{(1)}}\times\dots\times\Sym_{\lambda^{(\ell)}}$ be the
corresponding \textbf{Young subgroup}, or \textbf{parabolic subgroup},
of~$\Sym_n$. In particular, $\Row(\tlam)=\Sym_\blam$ and
$\Col(\tllam)=\Sym_{\blam'}$, where~$\Sym_\blam$ and~$\Sym_{\blam'}$
are the natural subgroups of~$\Sym_n$ associated with the
multipartitions $\blam$ and~$\blam'$.

For $\blam\in\Parts$ define $m_\blam=u^\blam x^\blam$, where \[ u^\blam
=\prod_{m=1}^n \prod_{l=\comp_m(\tlam)+1}^\ell(L_m-[\kappa_l])
\quad\text{and}\quad x^\blam=\sum_{w\in\Row(\tlam)}\xi^{\ell(w)}T_w.  \]
Let $M^\blam=m_\blam\Hn$ and set $M=\bigoplus_{\blam\in\Parts}M^\blam$.
Below we define the cyclotomic Schur algebra to be the endomorphism
algebra of~$M$.

Fix $\blam\in\Parts$ and define $m_{\s\t}=T_{d(\s)}^*m_\blam T_{d(\t)}$,
for $\s,\t\in\Std(\blam)$.  By \cite[Theorem 3.26]{DJM:cyc},
\[ \set{m_{\s\t}|\s,\t\in\Std(\blam)\text{ and }\blam\in\Parts} \]
is a cellular basis of $\Hn$, where $\Parts$ is ordered by dominance.
Consequently, if $\Hlam$ is the $K$-subspace of $\Hn$ spanned by
$\set{m_{\s\t}|\s,\t\in\Std(\bmu)\text{ for some }
               \bmu\in\Parts \text{ with }\bmu\gdom\blam}$,
then $\Hlam$ is a two-sided ideal of~$\Hn$.

Invoking the theory of cellular algebras~\cite{GL,M:Ulect}, for
$\blam\in\Parts$ there is a right cell module~$S^\blam$, called a
\textbf{Specht module}. By definition,
$S^\blam\cong m_{\tlam\tlam}\Hn/(m_{\tlam\tlam}\Hn\cap\Hlam)$.  Set
$m_\t=m_{\tlam\t}+\Hlam$, for $\t\in\Std(\blam)$. Then
$\set{m_\t|\t\in\Std(\blam))}$ is a $K$-basis for~$S^\blam$. The Specht
module comes equipped with an associative bilinear form
$\<\ ,\ \>\map{S^\blam\times S^\blam}K$ that is determined by
\begin{equation}\label{E:BilinearForm}
    m_\s m_{\t\u} = \<m_\s,m_\t\>m_\u,\qquad\text{ for all }
       \s,\t,\u\in\Std(\blam).
\end{equation}
Let $\rad S^\blam=\set{x\in S^\blam|\<x,y\>=0\text{ for all }y\in S^\blam}$.
Then $\rad S^\blam$ is an $\Hn$-submodule of~$S^\blam$. We define $D^\blam=S^\blam/\rad S^\blam$.

We need a parallel construction for the dual Specht modules. For $\blam\in\Parts$
set $n_\blam=u_\blam x_\blam$, where
\[ u_\blam =\prod_{m=1}^n \prod_{l=1}^{\comp_m(\tllam)-1}(L_m-[\kappa_l])
\quad\text{and}\quad
   x_\blam=\sum_{w\in\Col(\tllam)}(-\xi)^{-\ell(w)}T_w.
\]
For $\s,\t\in\Std(\blam)$ set $n_{\s\t}=T_{d'(\s)}^{\ast}n_\blam T_{d'(t)}$.
By \cite[(2.7)]{DuRui:branching},
$\set{n_{\s\t}|\s,\t\in\Std(\blam)\text{ and }\blam\in\Parts}$
is a second cellular basis of $\Hn$, where $\Parts$ is ordered
by \textit{reverse} dominance.

For each multipartition  $\blam\in\Parts$ the theory of cellular
algebras gives us a \textbf{dual Specht module}
$S_\blam=n_{\tllam\tllam}\Hn/(n_{\tllam\tllam}\Hn\cap\Hllam)$. Here, $\Hllam$ is
the two-sided ideal of~$\Hn$ with basis
\[\set{n_{\s\t}|\s,\t\in\Std(\bmu)\text{ for some }
               \bmu\in\Parts \text{ with }\blam\gdom\bmu},\]
Much as before, set $D_\blam=S_\blam/\rad S_\blam$.

\begin{Remark}\label{R:nbasislabels}
  The labelling that we are using for the $n$-basis is conjugate to
  that used in ~\cite{M:tilting,M:gendeg,DuRui:branching}
  so~$n_{\blam}$ should be replaced with $n_{\blam'}$, and $n_{\s\t}$ with
  $n_{\s'\t'}$, when comparing with these papers. Our notation
  reflects that fact that the elements $m_\blam$ and~$n_\blam$ come from
  looking at the row and column stabilizers of~$\tlam$ and~$\tllam$,
  respectively. Similarly, our labelling of the dual Specht
  modules follows the same convention, which is in agreement with the
  papers \cite{KMR:UniversalSpecht,HuMathas:QuiverSchurI}. In \cite{M:tilting}
  the module $S_\blam$ is written as $S_{\blam'}$.
\end{Remark}

Implicitly, the definitions of the elements $m_\blam$, $n_\blam$,
$m_{\s\t}$ and $n_{\s\t}$ all depend upon the choice of Hecke
parameter~$\xi$ and the multicharge $\charge$. This remark will be
important below when we vary these parameters.
For future use we summarise the properties of these modules that
follow directly  from the general theory of cellular algebras.

\begin{Theorem}[\cite{GL,DJM,M:tilting}]\label{T:CellularAlgebras}
  Suppose that $K$ is a field and $n\ge0$.
  \begin{enumerate}
    \item The $\set{D^\bmu|\bmu\in\Parts\text{ and }D^\bmu\ne0}$ is a complete
    set of pairwise non-isomorphic irreducible $\Hn$-modules.
    Moreover, $(D^\bmu)^*\cong D^\bmu$.
    \item The $\set{D_\bmu|\bmu\in\Parts\text{ and }D_\bmu\ne0}$ is a complete
    set of pairwise non-isomorphic irreducible $\Hn$-modules.
    Moreover, $(D_\bmu)^*\cong D_\bmu$.
  \end{enumerate}
\end{Theorem}

For the main results in this paper we need to describe the Specht modules
and dual Specht modules as submodules of $\Hn$, which is already known, and we
need to determine the isomorphisms between the two sets of simple
modules given by \autoref{T:CellularAlgebras}.

We extend the dominance ordering to the set of all standard
tableaux by defining $\s\gedom\t$ if
$$\Shape(\s_{\downarrow m})\gedom\Shape(\t_{\downarrow m}),$$
for $1\le m\le n$.

As remarked above, if $\ell>2$ then no pairs of dual bases for~$\Hn$ are
known. The following fundamental result implies that the two bases
$\set{m_{\s\t}}$ and $\set{n_{\s\t}}$ are dual bases of~$\Hn$ ``modulo
higher terms''.

\begin{Proposition}[%
  \protect{Hu and Mathas~\cite[Corollary~2.10]{HuMathas:GradedInduction}}, {Mathas~\cite[Theorem 5.5]{M:tilting}}]
  \label{P:mnDual}
  Suppose that $(\s,\t),(\u,v)\in\Std^2(\Parts)$ are pairs of tableaux
  of the same shape.  Then~$m_{\s\t}n_{\t\s}\ne0$ and
  $m_{\s\t}n_{\u\v}\ne0$ only if $\u\gedom\t$. Similarly, $n_{\t\s}m_{\s\t}\ne0$ and
  $n_{\u\v}m_{\s\t}\ne0$ only if $\v\gedom\s$.
\end{Proposition}

For $\blam\in\Parts$ set $z^\blam=n_\blam T_{w_{\blam'}}m_{\blam}$ and
$z_\blam=m_\blam T_{w_{\blam}}n_\blam$. Observe that $z^\blam=z_\blam^*$
since $w_{\blam}^{-1}=w_{\blam'}$. Moreover,
$z^\blam=n_{\tllam\tlam}m_{\tlam\tlam}=n_{\tllam\tllam}m_{\tllam\tlam}\ne0$ and
$z_\blam=m_{\tlam\tllam}n_{\tllam\tllam}=m_{\tlam\tlam}n_{\tlam\tllam}\ne0$
by \autoref{P:mnDual}.

\begin{Lemma}[\protect{Du and Rui~\cite[Theorem~2.9]{DuRui:branching}}]
  \label{L:SpechtSubmodules}
  Suppose that $\blam\in\Parts$. Then
      $S^\blam\cong z^\blam\Hn$ and
      $S_\blam\cong z_\blam\Hn$ as right $\Hn$-modules.
\end{Lemma}

To prove this observe that, because $z^\blam = n_\blam
T_{w_{\blam'}}m_{\tlam\tlam}$, there is a well-defined $\Hn$-module
homomorphism $\alpha\map{S^\blam}z^\blam\Hn$ given by
$\alpha(m_{\tlam}h)=z^\blam h$, for $h\in\Hn$. To complete the proof it remains to check that
$\set{z^\blam T_{d(\t)}|\t\in\Std(\blam)}$ is a basis of~$z^\blam\Hn$.
See \cite[Theorem~2.9]{DuRui:branching} or
\cite[Proposition~3.13]{M:gendeg} for details. The proof that
$z_\blam\Hn\cong S_\blam$ is similar.

As their names suggest, the Specht modules and dual Specht modules are
dual to each other. The proof of this requires the trace form~$\tau$
from \autoref{T:SymmetricAlgebra} and a strengthening of \autoref{P:mnDual}.

\begin{Theorem}[\protect{Mathas~\cite[Theorem~5.9]{M:gendeg}}]
  \label{T:NonvashingTrace}
  Suppose that $\blam\in\Parts$. Then
  $\tau(z^\blam T_{w_\blam})=\tau(z_\blam T_{w_{\blam'}})\ne0$
\end{Theorem}

We can now prove that $S^\blam$ and $S_\blam$ are dual to each other.

\begin{Corollary}[\protect{Mathas~\cite[Corollary~5.7]{M:tilting}}]
  \label{C:DualSpechts}
  Let $\blam\in\Parts$. Then $S^\blam\cong(S_\blam)^*$
  and $S_\blam\cong(S^\blam)^*$.
\end{Corollary}

\begin{proof}
  This result is proved in \cite{M:tilting} but we give a proof
  of the isomorphism  $S_\blam\cong(S^\blam)^*$ because we need the details below. For
  $\s\in\Std(\blam)$ let $\theta_\s\in (S^\blam)^*$ be the linear map
  determined by $\theta_\s(m_\t)=\tau(n_{\tllam\s}m_{\t\tlam}T_{w_{\blam}})$, for all
  $\t\in\Std(\blam)$. Now define a linear map $\theta\map{S_\blam}
  (S^\blam)^*$ by $\theta(n_\s)=\theta_\s$, for all $\s\in\Std(\blam)$.
  Then $\theta$ is a vector space isomorphism by \autoref{P:mnDual} and
  \autoref{T:NonvashingTrace}. For $h\in\Hn$ write $n_\s h=\sum_\v r_\v n_\v$,
  for $r_\v\in K$. Fix $\t\in\Std(\blam)$.  Using \autoref{P:mnDual},
  \begin{align*}
  \theta(n_\s h)(m_\t)
             &=\sum_{\v\in\Std(\blam)} r_\v\theta_\v(m_\t)
              =\sum_{\v\in\Std(\blam)} r_\v\tau(n_{\tllam\v}m_{\t\tlam}T_{w_{\blam}})\\
             &=\tau(n_{\tllam\s}hm_{\t\tlam}T_{w_{\blam}})
              =\theta_\s(m_\t h^*)
              =\bigl(\theta(n_\s) h\bigr)(m_\t).
  \end{align*}
   Hence, $\theta$ is an $\Hn$-module homomorphism, completing the proof.
\end{proof}

In \autoref{S:Tilting} below we need an analogue of \autoref{C:DualSpechts}
relating the simple modules $D^\bmu$ and $D_\bnu$, for
$\bmu,\bnu\in\Parts$.  To establish this we use the following
characterisation of the simple $\Hn$-modules together with the ``dual''
algebras~$\Hn'$ that are introduced below. This description of the simple
$\Hn$-modules as submodules of~$\Hn$ generalises a remark made by
James~\cite[p.~41]{James} for the symmetric groups.

\begin{Lemma}\label{L:SimpleSubmodules}
  Suppose that $\bmu\in\Parts$.
  \begin{enumerate}
    \item The simple module $D^\bmu\ne0$ if and only if $z^\bmu\Hn m_\bmu\ne0$.
    \item The simple module $D_\bmu\ne0$ if and only if $z_\bmu\Hn n_\bmu\ne0$.
  \end{enumerate}
  Moreover, if $D^\bmu\ne0$ then $D^\bmu\cong z_\bmu T_{w_{\bmu'}}m_\bmu\Hn$
  and if $D_\bmu\ne0$ then $D_\bmu=z^\bmu T_{w_{\bmu}}n_\bmu\Hn$.
\end{Lemma}

\begin{proof}We consider only the claims for $D^\bmu$ since those for
  $D_\bmu$ can be proved in the same way. As remarked after
  \autoref{L:SpechtSubmodules}, the map $S^\bmu\longrightarrow z^\bmu\Hn$
  determined by $m_\t\mapsto z^\bmu T_{d(\t)}$, for $\t\in\Std(\bmu)$,
  is an isomorphism. If $\s,\t\in\Std(\bmu)$ then using the
  definitions
  \[
     \<m_\s,m_\t\>z^\bmu
        = \<m_\s,m_\t\>n_{\t_\bmu\tmu} m_{\tmu\tmu}
        = n_{\t_\bmu\tmu} m_{\tmu\s}m_{\t\tmu}
        =z^\bmu T_{d(\s)}T_{d(\t)}^* m_\bmu,
  \]
  where the second equality follows from \autoref{E:BilinearForm} because
  $n_{\t_\bmu\tmu}\Hlam[\bmu]=0$ by \autoref{P:mnDual}. Hence, $D^\bmu\ne0$ if and
  only if $z^\bmu\Hn m_\bmu\ne0$, proving~(a).

  For the second claim, suppose that $D^\bmu\ne0$. We need to prove that
  $D^\bmu\cong z_\bmu T_{w_{\bmu'}}m_\bmu\Hn$.
  Combining \autoref{L:SpechtSubmodules} and \autoref{C:DualSpechts}, there
  are $\Hn$-module homomorphisms
  \[
       S^\bmu \bijection[\alpha] z^\bmu\Hn
              \bijection[\beta] z_\bmu\Hn
              \bijection[\theta](S^\bmu)^*,
  \]
  where $\alpha(m_\s)=z^\bmu T_{d(\s)}$, $\beta(x)=m_\bmu T_{w_\bmu}x$
  and $\theta(z_\bmu T_{d(\s)})=\theta_\s$, for $\s\in\Std(\bmu)$ and
  $x\in z_\bmu\Hn$, where $\theta_\s$ is defined in the proof of \autoref{C:DualSpechts}.
  Moreover,~$\alpha$ and $\theta$ are both isomorphisms. Let
  $\vartheta=\theta\circ\beta\circ\alpha$. We claim that, up to a
  non-zero scalar, $\vartheta$ is the map
  $S^\bmu\longrightarrow(S^\bmu)^*$ induced by the inner product~$\<\ ,\
  \>$ on~$S^\bmu$. To see this fix $\s\in\Std(\bmu)$ and write
  \[     n_\bmu T_{w_\bmu'} m_\bmu T_{d(\s)}=n_\bmu m_{\t_\bmu\s}
                 =\sum_{\u,\v}c_{\u\v}n_{\u\v},
              \qquad\text{for }c_{\u\v}\in K.
  \]
  As $\set{n_{\u\v}}$ is a cellular basis, with $\Parts$ ordered by
  reverse dominance, it follows that $c_{\u\v}\ne0$ only if
  $\bmu\gedom\Shape(\u)=\Shape(\v)$ with equality only
  if $\u=\t_\bmu$.  If $\v\in\Std(\bmu)$ set $c_\v=c_{\t_\bmu\v}$. Using
  \autoref{P:mnDual} for the third equality,
  \begin{align*}
      \vartheta(m_\s)&= \theta(z_\bmu T_{w_{\bmu'}} m_\bmu T_{d(\s)})
        =\theta\Bigl(\sum_{\u,\v} c_{\u\v} m_{\tmu\t_\bmu}n_{\u\v}\Bigr)
        =\theta\Bigl(\sum_{\v\in\Std(\bmu)} c_{\v} m_{\tmu\t_\bmu}n_{\t_\bmu\v}\Bigr)
        =\sum_{\v\in\Std(\bmu)} c_{\v}\theta_\v.
  \end{align*}
  Therefore, if $\t\in\Std(\bmu)$ then, using \autoref{P:mnDual} and
  \autoref{E:BilinearForm} again,
  \begin{align*}
    \vartheta(m_\s)(m_\t)
        &=\sum_{\v\in\Std(\bmu)}c_\v\theta_\v(m_\t)
         =\sum_{\v\in\Std(\bmu)}c_\v\tau(n_{\t_\bmu\v}m_{\t\tmu}T_{w_{\bmu}})
         =\tau\Bigl(\sum_{\u,\v}c_{\u\v}n_{\u\v}m_{\t\tmu}T_{w_{\bmu}}\Bigr)\\
         &=\tau\Bigl(n_\bmu T_{w_\bmu'} m_{\tmu\s}m_{\t\tmu}T_{w_{\bmu}}\Bigr)
          =\<m_\s,m_\t\>\tau\Bigl(n_\bmu T_{w_\bmu'} m_{\tmu\tmu}T_{w_{\bmu}}\Bigr)
          =\<m_\s,m_\t\>\tau(z^\bmu T_{w_{\bmu}}).
  \end{align*}
  By \autoref{T:NonvashingTrace}, $\tau(z^\bmu T_{w_{\bmu}})=\tau(z_\bmu
  T_{w_{\bmu'}})$ is a non-zero element of~$K$.  Hence, up to an
  invertible scalar, the map~$\vartheta$ coincides with the natural map
  induced by the bilinear form on~$S^\bmu$. In particular, if
  $D^\bmu\ne0$ then $D^\bmu\cong\im\vartheta\cong\im(\beta\circ\alpha)$.
  By definition, $\im(\beta\circ\alpha)=z_\bmu T_{w_\bmu'}m_\bmu\Hn$ so
  this completes the proof.
\end{proof}

To identify the isomorphic simple modules $D^\bmu$ and $D_\bnu$ we need
to shift to a ``dual'' Hecke algebra. As we make precise in
\autoref{T:RingelDuality} below, this is a shadow of Ringel duality in the
Hecke algebra world.

Let $\Hn'=\Hn^{\charge'}(\xi^{-1})$ be the cyclotomic Hecke algebra with
Hecke parameter $\xi^{-1}$ and multicharge
$\charge'=(-\kappa_\ell,\dots,-\kappa_1)$.  To help distinguish between
the elements of the algebras $\Hn$ and $\Hn'$ let
$T_1',\dots,T_{n-1}',L_1',\dots,L_n'$ be the generators of~$\Hn'$. More
generally, we decorate all elements of~$\Hn'$ with an appropriately
placed ${}'$. For example, $m_\blam'=u^\blam_{\prime}x^\blam_{\prime}$
and $n_\blam'=u'_\blam x'_\blam$ both belong to $\Hn'$. Similarly, there
are cellular bases $\set{m'_{\s\t}}$ and $\set{n'_{\s\t}}$ that give
rise to Specht modules $S^\blam_\prime$ and dual Specht modules~$S^\prime_\blam$,
for $\blam\in\Parts$. Let $D^\blam_\prime=S^\blam_\prime/\rad S^\blam_\prime$
and $D_\blam^\prime=S_\blam^\prime/\rad S_\blam^\prime$.

The next result, which does not appear to be in the literature,
generalises the sign automorphism on the group algebra of a symmetric
group. We leave the proof as an exercise to the reader because it
follows by simply inspecting the relations in the two algebras $\Hn$
and~$\Hn'$. (Recall from after \autoref{E:QInt} that
$[k]_\xi=-[-k]_{\xi^{-1}}$, for all $k\in\Z$.)

\begin{Lemma}\label{L:HashInvolution}
There is a unique algebra
isomorphism $\#:\Hn'\bijection\Hn$ such that
\[
    (T'_r)^\#=-T_r\qquad\text{and}\qquad
    (L'_s)^\#=-L_s, \qquad\quad\text{for }1\le r<n\text{ and }1\le s\le n.
\]
Moreover, $m_{\blam}=\pm(n'_{\blam'})^\#$,
$n_{\blam}=\pm(m'_{\blam'})^\#$,
$m_{\s\t}=\pm(n'_{\s'\t'})^\#$,
$n_{\s\t}=\pm(m'_{\s'\t'})^\#$, for all
$\blam\in\Parts$ and all $\s,\t\in\Std(\blam)$.
\end{Lemma}

The isomorphism $\#\map{\Hn'}\Hn$ induces an equivalence of
categories $\HEquiv\map{\RMod\Hn}{\RMod\Hn'}$. To describe the effect of
$\HEquiv$ on the simple $\Hn'$-modules define
\[\Klesh=\set{\bmu\in\Parts|D^\bmu\ne0}\quad\text{and}\quad
  \Klesh'=\set{\bmu\in\Parts|D_{\prime}^\bmu\ne0}.
\]
Recall from after \autoref{E:QInt} that $\xi$ and $\xi^{-1}$ both have
quantum characteristic~$e$. By the results in
\cite{Ariki:class,BK:GradedDecomp}, $\blam\in\Klesh$ if and only if
$\blam$ is a \textbf{Kleshchev}, or \textbf{restricted}, multipartition
with respect to~$(e,\charge)$ and $\blam\in\Klesh'$ if and only if
$\blam$ is Kleshchev with respect to~$(e,\charge')$. The Kleshchev
multipartitions index the crystal graphs of certain integrable highest
weight modules, so their definition is explicit but
recursive~\cite{Ariki:class}. The crystal graph combinatorics shows that
there is a crystal isomorphism, known as the \textbf{Mullineux map},
\[\bm\map{\Klesh}\Klesh'; \bmu\mapsto\bm(\bmu).\]
(Using the set $I$ defined in \autoref{S:FayersConjecture}, if $\bi\in I^n$
labels a path from $\bzero=(0|\dots|0)$ to $\bmu$ in the crystal graph
$\bigcup_{m\ge0}\Klesh[m]$ then $-\bi$ labels a path from $\bzero$
to~$\bm(\bmu)$ in the crystal graph $\bigcup_{m\ge0}\Klesh[m]'$. This
property determines the map $\bm$ uniquely. See
\cite[Theorem~4.12]{BK:GradedDecomp} or \cite{Ariki:branching} for more
details.) The next result connects the Mullineux map $\bm$ with
the representation theory of~$\Hn$.

\begin{Proposition}[cf.
  \cite{BK:GradedDecomp,KMR:UniversalSpecht,Mathas:Singapore}]
  \label{P:Mullineux}
  Suppose that $\blam\in\Parts$ and $\bmu\in\Klesh$. Then
  \[
     S_\prime^{\blam'}\cong\HEquiv(S_{\blam}),\quad
     S^\prime_{\blam'}\cong\HEquiv(S^{\blam})\quad\text{and}\quad
     D_\prime^{\bm(\bmu)}\cong \HEquiv(D^{\bmu})
  \]
  as right $\Hn'$-modules,
\end{Proposition}

\begin{proof}
  For graded Specht modules the first two isomorphisms are proved in
  \cite[Theorem~8.5]{KMR:UniversalSpecht}, although it is not completely
  clear that the isomorphism in~\autoref{L:HashInvolution} agrees with the
  homogeneous sign isomorphism considered in \cite{KMR:UniversalSpecht}.
  Fortunately, these two isomorphisms follow directly from the
  definitions because
  $m'_{\s'\t'}=\pm(n_{\s\t})^\#$ and $n'_{\s'\t'}=\pm(m_{\s\t})^\#$
  by \autoref{L:HashInvolution}. The isomorphism
  $S_\prime^{\blam'}\cong\HEquiv(S_{\blam})$, together with the modular
  branching rules~\cite{Ariki:branching,Groj:control} and a standard
  argument due to Kleshchev, now imply that
  $D_\prime^{\bm(\bmu)}\cong\HEquiv(D^\bmu)$. For the corresponding
  result for the graded simple modules see
  \cite[Theorem~3.6.6]{Mathas:Singapore} and \cite[(3.53)]{BK:GradedDecomp}.
\end{proof}

We can now prove the promised comparison results for the
simple $\Hn$-modules.

\begin{Corollary}\label{C:SimpleIsomorphism}
  Suppose that $\bmu\in\Parts$. Then $D_\bmu\ne0$ if and only if
  $\bmu'\in\Klesh'$. Moreover, if $\bmu\in\Klesh$ then
  $D^\bmu\cong D_{\bm(\bmu)'}$ as $\Hn$-modules.
\end{Corollary}

\begin{proof}
  By \autoref{L:SimpleSubmodules}(b) and \autoref{L:HashInvolution}, $D_\bmu\ne0$ if and only if
  $D^{\bmu'}_\prime\ne0$, so that $D_\bmu\ne0$ if and only if
  $\bmu'\in\Klesh'$. Hence, using \autoref{P:Mullineux} twice, if
  $\bmu\in\Klesh$ then
  \[ \HEquiv(D_{\bm(\bmu)'})\cong\HEquiv(\hd S_{\bm(\bmu)'})
                    \cong\hd\HEquiv(S_{\bm(\bmu)'})
                    \cong\hd S_\prime^{\bm(\bmu)}
                    \cong D_\prime^{\bm(\bmu)}
                    \cong\HEquiv( D^{\bmu}).
  \]
  Hence, $D_{\bm(\bmu)'}\cong D^\bmu$ as required.
\end{proof}

\begin{Corollary}\label{C:HDmu}
  Suppose that $\bmu\in\Parts$. Then $D_\bmu^\prime\ne0$ if and only if
  $\bmu'\in\Klesh$. Moreover, if  $\bmu'\in\Klesh$
  then $D^\prime_\bmu\cong D_\prime^{\bm(\bmu')}$ as $\Hn'$-modules.
\end{Corollary}

\begin{proof}
  By \autoref{L:SimpleSubmodules}(a) and \autoref{L:HashInvolution}, $D^\prime_\bmu\ne0$ if and
  only if $D^{\bmu'}\ne0$, so $D'_\bmu\ne0$ if and only if
  $\bmu'\in\Klesh$. Now suppose that $\bmu'\in\Klesh$. Then
  \[
      D^\prime_\bmu=\hd S^\prime_\bmu
                   \cong\hd\HEquiv(S^{\bmu'})
                   =\HEquiv(\hd S^{\bmu'})
                   \cong\HEquiv(D^{\bmu'})
                   \cong D_\prime^{\bm(\bmu')},
  \]
  where the second isomorphism follows by \autoref{P:Mullineux}.
\end{proof}

\begin{Corollary}\label{C:SimpleSocles}
  Suppose that $\bmu\in\Klesh$. Then
  $\soc S_{\bmu}\cong D^{\bmu}\cong \soc S^{\bm(\bmu)'}$.
\end{Corollary}

\begin{proof}
    If $\bmu\in\Klesh$ then $\hd S^\bmu\cong D^\bmu$. Moreover, $\hd S_{\bm(\bmu)'}\cong D_{\bm(\bmu)'}\cong D^\bmu$ by \autoref{C:SimpleIsomorphism}. Now take duals using \autoref{C:DualSpechts}
    (and \autoref{T:CellularAlgebras}).
\end{proof}

\section{Cyclotomic Schur algebras}\label{S:SchurAlgebras}
  We are now ready to introduce the cyclotomic Schur algebras, which are
  one of the main objects of study in this paper. Our first goal is to
  describe which simple modules appear in the socles of the Weyl
  modules, tilting modules and the projective indecomposable modules
  and, as an application describe the socles of the Young modules. In
  the last section we saw that the Specht modules and dual Specht
  modules came from two different cellular bases for~$\Hn$ and we needed
  to introduce the algebra $\Hn'$ to understand the difference in the
  labelling of the simple modules coming from these two bases.
  Similarly, in this section we introduce two variations of the
  cyclotomic Schur algebras which we will use to translate information
  about the socle of a tilting module to the socle of a projective
  module.

  As a prelude to defining the Schur algebras, for $\blam\in\Parts$ let
\begin{align*}
    M^\blam&=m_\blam\Hn,  &
    N^\blam&=n_\blam\Hn,  &
    M_{\prime}^\blam&=m'_\blam\Hn'  &
    N_{\prime}^\blam&=n'_\blam\Hn',
\intertext{and set}
     M&=\bigoplus_{\blam\in\Parts}M^\blam,  &
     N&=\bigoplus_{\blam\in\Parts}N^\blam,  &
     M_\prime&=\bigoplus_{\blam\in\Parts}M_\prime^\blam  &
     N_\prime&=\bigoplus_{\blam\in\Parts}N_\prime^\blam.
\end{align*}
Then $M^\blam$, $N^\blam$, $M$ and $N$ are $\Hn$-modules and
$M_\prime^\blam$, $N_\prime^\blam$, $M_\prime$ and $N_\prime$ are $\Hn'$-modules.

\begin{Proposition}\label{P:SchurPairs}
  Suppose that $n\ge0$. Then $(\Hn,M)$, $(\Hn,N)$, $(\Hn',M_{\prime})$ and
  $(\Hn',N_{\prime})$ are self-dual Schur pairs.
\end{Proposition}

\begin{proof}
  By \autoref{T:SymmetricAlgebra},  $\Hn$ and $\Hn'$ are both symmetric
  algebras and so, in particular, they are self-injective. Moreover,
  since $\Hn\cong M^{(0|\dots|0|1^n)}\cong N^{(n|0|\dots|0)}$ as right
  $\Hn$-modules, $M$ and $N$ are both faithful $\Hn$-modules. Therefore,
  $(\Hn,M)$ and $(\Hn,N)$ are both Schur pairs, Similarly, $(\Hn',M')$
  and $(\Hn',N')$ are both Schur pairs. It remains to show that these
  Schur pairs satisfy the conditions of \autoref{D:SelfDualSchurPair}.

  As algebras $\Hn$ and $\Hn'$ are both cellular. They have
  an anti-isomorphism, which is the unique
  anti-isomorphism that fixes each of the generators of~$\Hn$ and
  of~$\Hn'$, respectively. Moreover, by \autoref{T:CellularAlgebras}, the simple
  $\Hn$-modules and the simple $\Hn'$-modules are self-dual
  with respect to these involutions. To show that $(\Hn,M)$, $(\Hn,N)$,
  $(\Hn',M_\prime)$ and $(\Hn',N_\prime)$ are self-dual Schur pairs it is enough
  to show that each of the modules $M^\blam$, $N^\blam$,
  $M_\prime^\blam$ and $N_\prime^\blam$ is self-dual, for
  $\blam\in\Parts$. This is proved in \cite[Proposition~5.13]{M:tilting}.
  (The careful reader will notice that \cite[Proposition~5.13]{M:tilting} assumes
  that the cyclotomic parameters $Q_1,\dots,Q_\ell$ are
  invertible. As explained in \cite[\S1.1]{Mathas:Singapore}, because of
  \autoref{D:Hecke} this condition translates into the vacuous requirement
  that $[\kappa_r](\xi-\xi^{-1})+1=\xi^{2\kappa_r}$ is invertible
  in~$K$, for $1\le r\le\ell$. Hence, the requirement from
  \cite{M:tilting} that $Q_1,\dots,Q_\ell$ be invertible is satisfied.)
\end{proof}

Consequently, all of the results in \autoref{S:SchurPairs} apply to the
endomorphism algebras $S_A(M)$, where $(A, M)$ is one of the Schur pairs
given in \autoref{P:SchurPairs}. Notice that, by definition, $\Hn$ is a
direct summand of both $M$ and $N$ since
$M^{(0|\dots|0|1^n)}\cong\Hn\cong N^{(n|0|\dots|0)}$ as right
$\Hn$-modules. Similarly, $\Hn'$ is a direct summand of~$M_\prime$ and
of~$N_\prime$.

\begin{Definition}[Cyclotomic Schur
  algebras~\cite{DJM:cyc,BK:HigherSchurWeyl,M:tilting}]
  \label{D:Schur}
  The \textbf{cyclotomic Schur algebra} and the
  \textbf{twisted cyclotomic Schur algebra} of~$\Hn$ are the endomorphism
  algebras
  \[
      \Sn(M)=\End_{\Hn}(M)\qquad\text{and}\qquad
      \Sn'(M_\prime)=\End_{\Hn'}(M_\prime).
  \]
\end{Definition}

The cyclotomic Schur algebras include, as special cases, the
classical Schur algebras studied by Green~\cite{Green} and the $q$-Schur
algebras introduced by Dipper and James~\cite{DJ:Schur}. Just as our
definition of~$\Hn$ differs slightly from the definition given by Ariki
and Koike~\cite{AK}, the algebra $\Sn(M)$ is isomorphic to a cyclotomic
$\xi$-Schur algebra of \cite{DJM:cyc} when $\xi^2\ne1$ and it is isomorphic a
\textit{degenerate} cyclotomic Schur algebra~\cite{BK:HigherSchurWeyl}
when $\xi^2=1$.

We are mainly interested in the cyclotomic Schur algebra $\Sn=S_n(M)$,
however, we need the algebras $\Sn'=\Sn'(M_\prime)$ to prove some of our
results about $\Sn$-modules. Given \autoref{P:SchurPairs} it is also
natural to consider the algebras $\Sn(N)=\End_{\Hn}(N)$ and
$\Sn'(N_\prime)=\End_{\Hn'}(N_\prime)$. In fact, these algebras give no
additional information because \autoref{P:Mullineux} readily implies the
following.

\begin{Lemma}\label{L:SchurEquivs}
  If $\blam\in\Parts$ then $N^{\blam'}_\prime\cong\HEquiv(M^{\blam})$ and
  $M^{\blam'}_\prime\cong\HEquiv(N^{\blam})$ as right $\Hn$-modules.
  Therefore, $\Sn(M)\cong\Sn'(N_\prime)$ and $\Sn'(M_\prime)\cong\Sn(N)$.
\end{Lemma}

We leave the proof of \autoref{L:SchurEquivs} to the reader. In fact, as we
recall below, $\Sn=\Sn(M)$ and $\Sn'=\Sn'(M_\prime)$ are
quasi-hereditary algebras with weight posets~$\Parts$ and $\Parts^\op$,
respectively. Similarly, $\Sn(N)$ and $\Sn'(N_\prime)$ are
quasi-hereditary with weight posets $\Parts^\op$ and $\Parts$,
respectively. The isomorphisms $\Sn(M)\bijection\Sn'(N_\prime)$ and
$\Sn(N)\bijection\Sn'(M_\prime)$ of \autoref{L:SchurEquivs} are both
isomorphisms of quasi-hereditary algebras.

The algebras $\Sn$ and $\Sn'$ are finite dimensional cellular and
quasi-hereditary $K$-algebras with simple modules labelled by $\Parts$.
To describe cellular bases of~$\Sn$ and $\Sn'$ we need some more
combinatorics.

Fix $\bmu\in\Parts$ and $\blam\in\Parts$. A $\blam$-tableau of \textbf{type}
$\bmu$ is a map $\bT\map{[\blam]}\set{(k,r)|1\le
k\le\ell\text{ and } r\ge1}$ such that
$\mu_r^{(k)}=\#\set{\alpha\in[\blam]|\bT(\alpha)=(k,r)}$, for all $1\le
k\le\ell$ and $r\geq 1$. Order the pairs $\set{(k,r)}$
lexicographically.  Following~\cite[Definition~4.4]{DJM:cyc}, a
$\blam$-tableau~$\bT$ of type~$\bmu$ is \textbf{row semistandard} if
\begin{enumerate}
  \item the entries in $\bT$ are weakly increasing along rows,
  \item the entries in $\bT$ are strictly increasing down columns,
  \item if $\alpha=(k,r,c)\in[\blam]$ then $\bT(\alpha)\ge(k,r)$.
\end{enumerate}
Let $\RSStd(\blam,\bmu)$ be the set of row semistandard $\blam$-tableaux of
type $\bmu$. There is a bijection
$\Std(\blam)\bijection\RSStd(\blam,\bom)$ given by replacing each
entry~$m$ in a tableau with $(\ell,m)$, where $\bom=(0|\dots|0|1^n)$.

Let $\CSStd(\blam,\bmu)=\set{\bS'|\bS\in\RSStd(\blam',\bmu ')}$ be the
set of \textbf{column semistandard} $\blam$-tableaux of type~$\bmu$,
where conjugation reverses the order of components, swaps rows and
columns and replaces each entry~$(k,r)$ with $(\ell-k+1,r)$. (We
need these tableau in \autoref{S:Homomorphisms}.) If
$\bT\in\CSStd(\blam,\bmu)$ then the entries in~$\bT$ weakly increase
down columns, strictly increase along rows and if $(k,r,c)\in[\blam]$
then $\bT(\alpha)\le(k,c)$.  By
construction, if $\bT\in\CSStd(\blam,\bmu)$ then $\bT$ is a tableau of
type~$\bmu$.

Observe that $\RSStd(\blam,\bmu)\ne\emptyset$ only if $\blam\gedom\bmu$.
Similarly, $\CSStd(\blam,\bmu)\ne\emptyset$ only if
$\bmu\gedom\blam$. Let $\RSStd(\blam)=\bigcup_\bmu\RSStd(\blam,\bmu)$
and $\CSStd(\blam)=\bigcup_\bmu\CSStd(\blam,\bmu)$.

For $\t\in\Std(\blam)$ let $\row_\bmu(\t)$ be the tableau of type $\bmu$
obtained from $\t$ by replacing each entry~$m$ in~$\t$ with $(k,r)$ if
$m$ appears in \textit{row}~$r$ of the $k$th component of~$\tmu$ and define
$\col_\bmu(\t)$ similarly except that we use the \textit{column} index
of~$m$ in~$\t_\bmu$ instead.  It is not hard to see that every
semistandard tableau of type~$\bmu$ is equal to $\row_\bmu(\s)$ and to
$\col_\bmu(\t)$, for some $\s,\t\in\Std(\blam)$, although the converse is not
true in general. Let $\bT^{\blam}=\row_\blam(\tlam)$ and
$\bT_\blam=\col_\blam(\t_\blam)$, respectively, be the unique row and
column semistandard $\blam$-tableau of type~$\blam$. As a quick exercise
in these definitions, $\bT^\blam=\bT_{\blam'}'$.

By \cite[Theorem 4.14]{DJM:cyc}, $M^\bmu=m_{\bmu}\Hn$ has $K$-basis
$\set{m_{\bS\t}|\bS\in\RSStd(\blam,\bmu),\t\in\Std(\blam),\blam\in\Parts}$,
where
\[
m_{\bS\t}=\sum_{\substack{\s\in\Std(\blam)\\\row_\bmu(\s)=\bS}}\xi^{\ell(d(\s))}m_{\s\t}.
\]
Extending this idea, following~\cite{DJM:cyc}, if $\bS\in\RSStd(\blam,\bmu)$ and
$\bT\in\RSStd(\blam,\bnu)$ define
\begin{equation}\label{E:mST}
  m_{\bS\bT} = \sum_{\substack{\t\in\Std(\blam)\\\row_\bnu(\t)=\bT}}\xi^{\ell(d(\t))}m_{\bS\t}.
\end{equation}
Then $m_{\bS\bT}\in m_\bmu\Hn\cap\Hn m_\nu$ by the remarks above.
Therefore, the map $\phiST\map{M^\bnu}M^\bmu$ defined by
$\phiST(m_\bnu h)=m_{\bS\bT}h$, for all $h\in\Hn$, is an $\Hn$-module
homomorphism. We consider $\phiST$ as an element of~$\Sn$. By
\cite[Theorem 6.6]{DJM:cyc},
$\set{\phiST|\bS\in\RSStd(\blam,\bmu)\text{ and } \bT\in\RSStd(\blam,\bnu)
                \text{ for }\blam,\bmu,\bnu\in\Parts}$
is a cellular basis of~$\Sn$.

For each $\blam\in\Parts$ the algebra $\Sn$ has a \textbf{Weyl module}
$W^{\blam}$, which is the corresponding right cell module of~$\Sn$.
To make this more explicit, let $\Slam$ be the two-sided ideal
of~$\Sn$ with basis the $\phiST$ where $\bS$ and $\bT$ are
semistandard $\bmu$-tableau with $\bmu\gdom\blam$. Then
$W^{\blam}\cong\varphi_{\bT^{\blam}\bT^{\blam}} \Sn/%
        (\varphi_{\bT^{\blam}\bT^{\blam}} \Sn\cap\Slam)$.
If $\bS\in\RSStd(\blam)$ then set
$\varphi_\bS=\varphi_{\bT^{\blam}\bS}+\Slam\in W^\blam$. Then
$\set{\varphi_\bS|\bS\in\RSStd(\blam)}$
is a $K$-basis of $W^\blam$. Exactly as for a Specht module, the
Weyl module $W^\blam$ comes equipped with a bilinear form,
also written as $\<\ ,\ \>$, determined by
  \[\varphi_\bS \varphi_{\bU\bV} = \<\varphi_\bS,\varphi_\bU\>\varphi_\bV,
       \qquad\text{ for all } \bS,\bU,\bV\in\RSStd(\blam).\]
Let $L^\blam=W^\blam/\rad W^\blam$ where $\rad W^\blam$ is
the radical of the form~$\<\ ,\ \>$. By definition, $\varphi_{\Tlam\Tlam}$ is the
identity map on~$M^\blam$, so
$\<\varphi_{\Tlam},\varphi_{\Tlam}\>=1$ and, consequently, $L^\blam\ne0$.

In a similar way we can define elements $m'_{\bS\bT}\in M^\blam_\prime$
and homomorphisms
$\varphi_{\bS\bT}'\in\Hom_{\Hn'}(M^\bnu_\prime,M^\bmu_\prime)$, for
semistandard tableaux $\bS,\bT\in\RSStd(\blam)$. Hence, we can define
Weyl modules $W^\blam_\prime$ and simple modules
$L^\blam_\prime=W^\blam_\prime/\rad W^\blam_\prime$ for $\Sn'$.

The cellular algebra involutions on $\Hn$ and $\Hn'$ induce involutions $\ast$
on $\Sn$ and $\Sn'$ with the property that
$\varphi_{\bS\bT}\mapsto\varphi_{\bT\bS}$ and
$\varphi'_{\bS\bT}\mapsto\varphi'_{\bT\bS}$. For $\blam\in\Parts$ let
$V^\blam=(W^\blam)^*$ and
$V^\blam_\prime=(W^\blam_\prime)^*$ be \textbf{costandard
modules} for $\Sn$ and $\Sn'$, respectively.

\begin{Theorem}[\cite{DJM:cyc,M:tilting}]\label{T:CellularSAlgebras}
  Suppose that $n\ge0$.
  \begin{enumerate}
    \item The algebra $\Sn$ is a quasi-hereditary cellular algebra with
    weight poset $\Parts$,
    cellular basis
    \[\set{\varphi_{\bS\bT}|\bS\in\RSStd(\blam,\bsig)\text{ and }
                       \bT\in\RSStd(\blam,\btau)
                       \text{ for }\blam,\bsig,\btau\in\Parts},\]
    standard modules $\set{W^\blam|\blam\in\Parts}$,
    costandard modules $\set{V^\blam|\blam\in\Parts}$ and
    pairwise non-isomorphic self-dual simple modules $\set{L^\blam|\blam\in\Parts}$.
    \item The algebra $\Sn'$ is a quasi-hereditary cellular algebra with
    weight poset $\Parts^\op$,
    cellular basis
    \[\set{\varphi'_{\bS\bT}|\bS\in\RSStd(\blam,\bsig)\text{ and }
                       \bT\in\RSStd(\blam,\btau)
                       \text{ for }\blam,\bsig,\btau\in\Parts},\]
    standard modules $\set{W^\blam_\prime|\blam\in\Parts}$,
    costandard modules $\set{V^\blam_\prime|\blam\in\Parts}$ and
    pairwise non-isomorphic self-dual simple modules $\set{L^\blam_\prime|\blam\in\Parts}$.
  \end{enumerate}
\end{Theorem}

Although we will not need them, there are analogous definitions of
elements $n_{\bS\bT}\in N^\blam$ and $n_{\bS\bT}'\in N^\blam_\prime$ and maps
$\psi_{\bS\bT}\in\Hom_{\Hn}(N^\bmu,N^\bnu)$ and
$\psi_{\bS\bT}'\in\Hom_{\Hn'}(N^\bmu_\prime,N^\bnu_\prime)$, where
$\bS$ and $\bT$ are \textit{column} semistandard tableaux. Then
\autoref{L:SchurEquivs} and\autoref{T:CellularSAlgebras} imply that
$\set{\psi_{\bS\bT}}$ and $\set{\psi'_{\bS\bT}}$ are cellular bases for
the quasi-hereditary algebras $\Sn(N)$ and $\Sn'(N_\prime)$,
respectively.

Let $\pi$ be the natural projection from $M$ onto $\Hn$. Then it is easy to check that $\pi^\ast=\pi$. Recall that $(\Hn,M)$ is a Schur pair.
By the results in \autoref{S:SchurPairs}, we have that
\[
   S_n\pi\cong M=\bigoplus_{\blam\in\P_n}m_\blam\Hn,\quad
   \pi S_n\cong \bigoplus_{\blam\in\P_n}\Hn m_\blam\quad\text{and}\quad
  S_n^{\op}\cong\End_{\Hn}\Bigl(\bigoplus_{\blam\in\P_n}\Hn m_\blam\Bigr) .
\]
In particular, we can regard $\oplus_{\blam\in\P_n}\Hn m_\blam$ as a right $S_n$-module. Similar results hold for the Schur pair $(\Hn', M')$. In what follows we consider $z^\blam$ to be an element
of~$\Hn m_\blam$ and $z^\blam_{\prime}$ as an element of $\Hn'm'_\blam$.

The next result shows how the classical definitions of Weyl modules
extend to the cyclotomic case. In the semisimple case this follows
directly from \autoref{L:Ideals} and \autoref{L:SpechtSubmodules} but in
general we need to work harder. In the special case when $\ell=1$ this
is due to Dipper and James~\cite{DJ:qWeyl}.

\begin{Proposition}\label{P:WeylSubmodules}
Suppose that $\blam\in\Parts$.
\begin{enumerate}
  \item As $\Sn$-modules, $W^\blam\cong z^\blam\Sn $. In
  particular, $W^\blam$ is (isomorphic to) a submodule of~$\oplus_{\blam\in\P_n}\Hn m_\blam$.
  \item As $\Sn'$-modules, $W^\blam_\prime\cong z^\blam_{\prime}\Sn'$ . In
  particular, $W^\blam_\prime$ is (isomorphic to) a submodule of~$\oplus_{\blam\in\P_n}\Hn' m'_\blam$.
\end{enumerate}
\end{Proposition}

\begin{proof}
  We prove (a) and leave (b) for the reader. Recall that
  $W^\blam\cong\varphi_{\Tlam\Tlam}\Sn/(\varphi_{\Tlam\Tlam}\Sn\cap\Slam)$.
  Define an $\Sn$-module homomorphism by
  \[\theta\map{\varphi_{\Tlam\Tlam}\Sn}\bigoplus_{\blam\in\P_n}\Hn m_\blam;\,\,
          \phi\mapsto z^\blam\phi, \qquad\text{for all }\phi\in\varphi_{\Tlam\Tlam}\Sn.
  \]
  Fix semistandard tableaux $\bS\in\RSStd(\bmu,\blam),\bT\in\RSStd(\bmu)$, for some $\bmu\in\Parts$.
  Then
  \[\theta(\varphi_{\bS\bT})=z^\blam\cdot \phiST =(n_\blam T_{w_{\blam'}}m_{\blam})\phiST
                   =n_{\blam}T_{w_{\blam'}}m_{\bS\bT}.\]
  By \autoref{P:mnDual}, if $\s,\t\in\Std(\bmu)$ then
  $n_{\blam}m_{\s\t}\ne0$ only if $\blam\gedom\bmu$. Consequently,
  $\theta(\varphi_{\bS\bT})=0$ if $\blam\not\gedom\bmu$. In particular,
  $\Slam\cap\varphi_{\Tlam\Tlam}\Sn\subseteq\ker\theta$ since
  $\set{\phiST|\bS,\bT\in\RSStd(\bmu)\text{ where }\bmu\gdom\blam}$ is a basis
  of~$\Slam$. On the other hand, if $\bS\in\RSStd(\blam)$ then
  $\theta(\phiST[\Tlam\bS])=n_{\blam}T_{w_{\blam'}}m_{\tlam\bS}
                           =n_{\blam}m_{\tllam\bS}$.
  In view of \autoref{L:SpechtSubmodules},
  \[\set{\theta(\phiST[\Tlam\bS])|\bT\in\RSStd(\blam)}
            =\set{n_{\blam}m_{\tllam\bS}|\bS\in\RSStd(\blam)}\]
  is linearly independent (see also \cite[(2.10)]{DuRui:branching} or
  \cite[Proposition~5.9]{M:tilting}). Therefore,
  $W^\blam=\phiST[\Tlam\Tlam]\Sn/\ker\theta\cong\im\theta=z^\blam \Sn$
  as required.
\end{proof}

The last result in this section explains the significance of the algebras
$\Sn'\cong \Sn(N)$ in the representation theory of~$\Sn$ (and hence why
we need them in this paper). To do this we first recall the definition of
the Ringel dual of a quasi-hereditary algebra.

Let $S$ be a quasi-hereditary algebra with standard modules $W^i$
and costandard modules $V^i$, where~$i$ runs over a poset~$(I,\ge)$.
Let $\rDelMod{S}$ be the full subcategory of $\RMod{S}$ consisting of $\Delta$-filtered
$S$-modules. Thus, $X\in\rDelMod{S}$ if and only if $X$ has a filtration
with each subquotient isomorphic to a Weyl module $W^i$, for
$i\in I$. Similarly,  let $\rNabMod{S}$ be the full subcategory of $\RMod{S}$ consisting of
$\nabla$-filtered $S$-modules. If $X\in\rDelMod{S}$ let
$(X:W^i)$ be the number of subquotients of~$X$ isomorphic
to~$W^i$. Define $(Y:V^i)$ in the same way when
$Y\in\rNabMod{S}$. Since $S$ is quasi-hereditary the multiplicities
$(X:W^i)$ and $(Y:V^i)$ are independent of the choices
of $\Delta$ and $\nabla$ filtrations.

An \textbf{$S$-tilting module} is any $S$-module in
$\rDelMod{S}\cap\rNabMod{S}$. As $S$ is quasi-hereditary, by
\cite[A4]{Donkin:book} for each $i\in I$ there is a unique
indecomposable tilting module $T^i$ for~$S$ such that $(T^i:W^i)=1$
and $(T^i:W^j)\ne0$ only if $i\ge j$.  Moreover, up to isomorphism
$\set{T^i|i\in I}$ is a complete set of pairwise non-isomorphic
indecomposable tilting modules.

A \textbf{full tilting module} is a tilting module that contains
$\bigoplus_{i\in I} T^i$ as a summand. A \textbf{Ringel dual} of the
algebra~$S$ is any algebra $\RingelDual{S}=\End_{S}(T)$, where $T$ is any full
tilting module. Then $\RingelDual{S}$ is quasi-hereditary with respect to the opposite
poset~$I^\op$. By construction, the Ringel dual is unique up to Morita
equivalence. There is an exact functor $\Run\map{\rNabMod{S}}\rDelMod{\RingelDual{S}}$
given by $X\mapsto\Hom_S(T,X)$.

Returning now to the cyclotomic Schur algebras, let $T^\blam$ and
$T^\blam_\prime$ be the tilting modules for~$\Sn$ and~$\Sn'$, for
$\blam\in\Parts$. Let $P^\blam$ and $P^\blam_\prime$ be the projective
covers of~$L^\blam$ and~$L_\prime^\blam$, respectively.

Let $\Funn\map{\RMod\Sn}{\RMod\Hn}$ and $\Funn'\map{\RMod\Sn'}{\RMod\Hn'}$
be the Schur functors, defined in \autoref{S:SchurPairs}.

\begin{Theorem}[Ringel duality for cyclotomic Schur algebras]
  \label{T:RingelDuality}
  The twisted cyclotomic Schur algebra $\Sn'\cong\RingelDual{\Sn}$ is Ringel dual
  to~$\Sn$. Moreover, there is an exact functor $\Run\map{\rNabMod{\Sn}}\rDelMod{\Sn'}$
  such that the following diagram of functors commutes:
  \begin{center}
    \begin{tikzpicture}[>=stealth,->,shorten >=2pt,looseness=.5,auto]
      \matrix (M)[matrix of math nodes,row sep=1cm,column sep=16mm]{
          \rNabMod{\Sn} & \rDelMod{\Sn'}\\
          \RMod\Hn & \RMod\Hn'\quad.\\
       };
       \draw(M-1-1)--node[above]{$\Run$}(M-1-2);
       \draw(M-2-1)--node[below]{$\HEquiv$}(M-2-2);
       \draw(M-1-1)--node[left]{$\Funn$}(M-2-1);
       \draw(M-1-2)--node[right]{$\Funn'$}(M-2-2);
    \end{tikzpicture}
  \end{center}
  The functor $\Run$ is determined by
  $\Run(V^{\blam})\cong W_\prime^{\blam'}$, for all
  $\blam\in\Parts$. Moreover, $\Run(T^{\blam})\cong P_\prime^{\blam'}$
  as $\Sn$-modules.
\end{Theorem}

\begin{proof}
   By \cite[Corollary~7.3]{M:tilting}, there is a full tilting
   module~$E$ for $\Sn$ such that $\Sn(N)\cong\End_{\Sn}(E)$. By
   \autoref{L:SchurEquivs}, $\Sn'\cong\Sn(N)$ so $\Sn'$ is a Ringel dual
   of~$\Sn$.
   The isomorphism $\Sn'\bijection\Sn(N)$ induces an
   equivalence of categories $\RMod\Sn(N)\bijection\RMod\Sn'$, so Ringel
   duality give an exact functor $\Run\map{\rNabMod{\Sn}}{\rDelMod{\Sn'}}$
   that sends indecomposable tilting modules to indecomposable
   projective modules.
   By induction on the dominance order it follows that
   $\Run(V^{\blam})\cong W^{\blam'}_\prime$, for all
   $\blam\in\Parts$. Since $\Run$ is exact on $\rNabMod{\Sn'}$, it is
   uniquely determined by its action on the costandard modules. In
   particular, $\Run(T^{\blam})\cong P_\prime^{\blam'}$ since
   $V^{\blam}$ is a quotient of $T^\blam$.

   It remains to check that the diagram above commutes. By exactness it
   is sufficient to check commutativity on the costandard
   $\Sn$-modules. By \autoref{C:SchurFunctor}, \autoref{L:SpechtSubmodules}
   and \autoref{P:WeylSubmodules}, if $\blam\in\Parts$ then
   $\Funn(V^\blam)\cong S_\blam$ and
   $\Funn'(W_\prime^\blam)\cong S_\prime^\blam$. Therefore,
   as $\Hn$-modules,
   \[
   (\Funn'\circ\Run)(V^{\blam})
            \cong\Funn'(W_\prime^{\blam'})
            \cong S_\prime^{\blam'}
            \cong\HEquiv(S_{\blam})
            \cong (\HEquiv\circ\Funn)(V^{\blam}),
    \]
    where the third isomorphism comes from \autoref{P:Mullineux}.
\end{proof}

It is worth mentioning that the Ringel duality functor~$\Run$, when
considered a functor from $\RMod\Sn$ to $\RMod\Sn'$, is only
\textit{left\/} exact, and not right exact.  We restrict to the
subcategory $\rNabMod\Sn$ in \autoref{T:RingelDuality} only because this
ensures that $\Run$ is exact.

Following \cite[\S3-4]{M:tilting}, and \autoref{D:YoungModules}, we make
the following definition.

\begin{Definition}[\protect{\cite[\S3, \S4, \S7]{M:tilting}}] The
  \textbf{Young modules} and the \textbf{twisted Young modules} are the
  $\Hn$-modules $Y^\blam=\Funn(P^\blam)$ and $Y_\blam=\Funn(T^\blam)$,
  respectively, for $\blam\in\Parts$.  \end{Definition}

Similarly, define Young modules for~$\Hn'$ by setting
$Y^\blam_\prime=\Funn'(P^\blam_\prime)$ and
$Y_\blam^\prime=\Funn'(T^\blam_\prime)$. We need Ringel duality for the
next corollary, which we will use in \autoref{S:Tilting} to prove results
about the socles of Young modules for $\Hn$ and projective modules
for~$\Sn$.

\begin{Corollary}\label{C:EquivYoungModules}
  Suppose that $\blam\in\Parts$. Then
  $\HEquiv(Y_{\blam})\cong Y_\prime^{\blam'}$ as $\Hn'$-modules.
\end{Corollary}

\begin{proof}
   Using the definitions and \autoref{T:RingelDuality}, as $\Hn'$-modules,
   \[ \HEquiv(Y_{\blam})
            \cong(\HEquiv\circ\Funn)(T^{\blam})
            \cong(\Funn'\circ\Run)(T^{\blam})
            \cong\Funn'(P_\prime^{\blam'})
            \cong Y_\prime^{\blam'},
    \]
    as required.
\end{proof}

Let $I^\blam$ be the injective hull of $L^\blam$, for $\blam\in\Parts$.
Then $I^\blam\in\rNabMod{\Sn}$.

\begin{Corollary}\label{C:}
  Suppose that $\blam\in\Parts$. Then
  $\HEquiv(Y^{\blam})\cong Y^\prime_{\blam'}$.
\end{Corollary}

\begin{proof}
   The Schur functor commutes on the dualities on $\RMod\Sn$ and
   $\RMod\Hn$ and $(Y^\bmu)^*\cong Y^\bmu$ by
   \cite[Corollary~5.14]{M:tilting} (when $\blam\in\Klesh$ this follows
   from \autoref{L:YSelfDual}). Moreover, Ringel duality sends injective
   modules to tilting modules. Therefore,
   \[
       \HEquiv(Y^\blam)
            \cong\HEquiv\((Y^\blam)^*\)
            \cong(\HEquiv\circ\Funn)\((P^\blam)^*)
            \cong(\HEquiv\circ\Funn)(I^\blam)
            \cong\Funn'(T_\prime^{\blam'})
            \cong Y^\prime_{\blam'}.
    \]
\end{proof}

The final result in this section is part of the folklore for $\Hn$ but
as far as we are aware the result is not in the literature.

\begin{Corollary}
  Suppose that $\bmu\in\Parts$. Then:
  \begin{enumerate}
    \item If $\bmu\in\Klesh$ then $Y^\bmu$ is the projective cover of~$D^\bmu$.
    \item If $\bmu'\in\Klesh'$ then $Y_\bmu$ is the projective cover
          of~$D_{\bmu}$.
    \item If $\bmu\in\Klesh$ then $Y^\bmu\cong Y_{\bm(\bmu)'}$ as
  $\Hn$-modules.
\end{enumerate}
\end{Corollary}

\begin{proof}
  By \autoref{P:YoungModules}, $Y^\bmu$ is projective if and only if
  $\bmu\in\Klesh$, in which case it is the projective cover of~$D^\bmu$.
  This proves~(a). Similarly, $Y^{\bmu'}_\prime$ is projective if and only
  if $\bmu'\in\Klesh'$, in which case it is the projective cover
  of~$D^{\bmu'}_\prime\cong\HEquiv(D^{\bm^{-1}(\bmu')})$. Hence,~(b) follows
  by \autoref{C:EquivYoungModules} and \autoref{C:SimpleIsomorphism}. Part~(c) is
  now automatic from parts~(a) and~(b).
\end{proof}

\section{Socles of Weyl modules, tilting modules and projective modules}
\label{S:Tilting}

We now come to the first main result of this paper, which motivated
much of the development of \autoref{S:SchurPairs}. In more detail, we
describe the simple $\Sn$-modules that can appear in the socles of the
Weyl modules, tilting modules and projective indecomposable modules and
give similar results for the socles of the (twisted) Young modules.

Before we begin we note the following immediate consequence of
\autoref{P:SchurPairs} and \autoref{T:Equivalence}.

\begin{Theorem}\label{T:SnEquivalence}
  Suppose that  $\lambda\in\P_n$. Then
  the following are equivalent:
  \begin{enumerate}
    \item $\lambda\in\K_n$,
    \item $D^\lambda\ne0$,
    \item $L^\lambda$ is a right submodule of $M_R$,
    \item $P^\lambda$ is a direct summand of $M_R$,
    \item $P^\lambda$ is a projective-injective right $S_n$-module,
    \item $P^\lambda$ is an indecomposable tilting module,
    \item $P^\lambda$ is self-dual.
  \end{enumerate}
\end{Theorem}

There are analogous results for $\Sn'$-modules.

The first of these results generalises a classical result of
James~\cite[Theorem~2.8]{James:TensorDecomp} in level one (when
$\xi^2=1$ and $\ell=1$). This result will be used to prove Fayers'
conjecture in \autoref{T:Fayers} below.

\begin{Theorem} \label{T:DeltaSocle} Suppose that $\blam,\bmu\in\Parts$.
  \begin{enumerate}
    \item The simple module $L^{\bmu}$ is a submodule of~$W^{\blam}$ only
          if $\bmu\in\Klesh$.
    \item The simple module $L^{\bmu}_\prime$ is a submodule
    of~$W_\prime^{\blam}$ only if $\bmu\in\Klesh'$.
  \end{enumerate}
\end{Theorem}

\begin{proof}
  Both parts can be proved in the same way so we consider only~(a).  By
  \autoref{P:SchurPairs}, $(\Hn, M)$ is a Schur pair. Therefore, by
  \autoref{T:Socle}, $L^\bmu$ is an $\Sn$-submodule of~$M_R$ if and only if
  $D^\bmu\ne0$, which is if and only if $\bmu\in\Klesh$. On the other
  hand, $W^\blam$ is isomorphic to an $\Sn$-submodule of~$M_R$ by
  \autoref{P:WeylSubmodules}, so $\soc W^\blam\subseteq\soc M_R$.  Hence,
  $L^\bmu$ is a submodule of~$W^\blam$ only if $\bmu\in\Klesh$ as
  claimed.
\end{proof}

\begin{Corollary} \label{C:PTSocle} Suppose that $X\in\rDelMod\Sn$ and
  $\blam,\bmu\in\Parts$. Then $L^\bmu$ is a submodule of~$X$ only if
  $\bmu\in\Klesh$. In particular, $L^\bmu$ is a submodule
  of~$P^\blam\oplus T^\blam$ only if $\bmu\in\Klesh$.
\end{Corollary}

By \autoref{P:WeylSubmodules} and \autoref{L:Ideals}, or working directly
with the definitions, $S^\blam\cong\Funn(W^\blam)$ for all
$\blam\in\Parts$. Hence, applying the Schur functor to
\autoref{T:DeltaSocle} gives the following.

\begin{Corollary} \label{C:soc5} Suppose that $\blam,\bmu\in\Parts$. Then
  \[\relax
      [\soc W^\blam:L^\bmu]=[\soc S^\blam:D^\bmu]
        \qquad\text{and}\qquad
      [\soc T^\blam:L^\bmu]=[\soc Y_\blam:D^\bmu].
  \]
  In particular, these two multiplicities are non-zero only if $\bmu\in\Klesh$.
\end{Corollary}

\begin{proof} Recall from \autoref{P:SchurPairs} that $(\Hn,M)$ is a
  self-dual Schur pair.  Therefore, the result follows by
  \autoref{C:PTSocle} and \autoref{C:KleshSocle}.
\end{proof}

\begin{Corollary}
  Assume $\bmu\in\Klesh$. Then
  $\soc W^{\bm(\bmu)'}\cong L^{\bmu}$. Equivalently,
  $\hd V^{m(\bmu)'}\cong L^{\bmu}$.
\end{Corollary}

\begin{proof}
  Combine \autoref{C:soc5} and \autoref{C:SimpleSocles} for the first
  isomorphism and take duals for the second.
\end{proof}

\section{Proof of Fayers' conjecture}\label{S:FayersConjecture}

In this section, as our second application of the results in
\autoref{S:SchurPairs}, we prove Fayers Conjecture~\cite{Fayers:LLT}. This
is one of the main result of this paper. We first give a precise
statement of Fayers' conjecture and then recall some recent results that we
need from the graded representation theory of~$\Hn$ and~$\Sn$.

Throughout this section we assume that $K=\C$ and we fix $\xi\in\C$, an
element of quantum characteristic~$e$ (that is, a primitive $2e$th root
of unity in~$\C$), and a multicharge $\charge\in\Z^\ell$. In fact, we
will work with $\charge'=(-\kappa_\ell,\dots,-\kappa_1)$ because our
argument uses results of Stroppel and
Webster~\cite{StroppelWebster:QuiverSchur} who worked with the twisted
cyclotomic Schur algebra~$\Sn'\cong\Sn(N)$, by \autoref{L:SchurEquivs}.

Set $I=\Z/e\Z$, where we adopt the convention that
$e\Z=e\Z\cap\Z=\set{0}$ when $e=\infty$ (so that $I=\Z$ when
$e=\infty$).  Let $q$ be an indeterminate over $\Q$. Let $\Usl$ be the
\textbf{quantised enveloping algebra} of the \textbf{Kac-Moody algebra}
$\sl$ (in particular, we consider $\Usl[\infty]$ when $e=\infty$). The
algebra $\Usl$ is a Hopf algebra with coproduct $\Delta$ determined by
\[  \Delta(E_i)=E_i\otimes K_i+1\otimes E_i,\quad
    \Delta(E_i)=F_i\otimes1+ K_i^{-1}\otimes F_i,\quad\text{and}\quad
    \Delta(K_i)= K_i\otimes K_i,
\]
for $i\in I$. Let $\set{\Lambda_i|i\in I}$ be the set of
\textbf{fundamental weights} and $\set{\alpha_i|i\in I}$ the
\textbf{simple roots} for~$\sl$ and set $P^+=\bigoplus_{i\in
I}\N\Lambda_i$ and $Q^+=\bigoplus_{i\in I}\N\alpha_i$ Then $\Usl$ is the
$\Q(q)$-algebra generated by elements $\set{E_i, F_i, K_i^{\pm1}|i\in I}$
subject to the well-known quantised relations~\cite{Lusztig:QuantBook}. The
\textbf{bar involution}~$\BarInvolution$ is the $\Q$-linear ring
automorphism of ~$\Usl$ determined by
\[
   \overline{E}_i=E_i,\quad
   \overline{F}_i=F_i,\quad
   \overline{K}_i=K_i^{-1}\quad\text{and}\quad
   \overline{q}=q^{-1},
\]
for $i\in I$.

The (combinatorial) \textbf{Fock space} $\Fock$ is the $\Q(q)$-vector
space
\[\Fock=\bigoplus_{n\ge0}\bigoplus_{\blam\in\Parts}\Q(q)s_\blam.\]
The \textbf{residue} of a node $A=(k,r,c)\in[\blam]$
is $\res A=-\kappa_{\ell+1-k}+c-r+e\Z\in I$. If $\res A=i\in I$ then $A$ is an
\textbf{$i$-node}. Following
Hayashi~\cite{Hay} and Misra and Miwa~\cite{MisraMiwa}, the action of
$\Usl$ on $\Fock$ can be described explicitly using the combinatorics of
addable and removable $i$-nodes. At the categorical level the
$\Usl$-action corresponds to graded $i$-induction and $i$-restriction
for the cyclotomic Hecke
algebras~\cite{BKW:GradedSpecht,HuMathas:GradedInduction,Mathas:Singapore,
               BK:GradedDecomp}.
As we do not need the precise details we refer interested reader
to~\cite[\S3.6]{BK:GradedDecomp} or \cite[\S3.5]{Mathas:Singapore}.

Let $\Lambda'=\sum_{i\in I}l_i\Lambda_i$, where
$l_i=\#\set{1\le l\le\ell|i=-\kappa_l+e\Z}$, and set
$s_{\Lambda'}=s_{(0|\dots|0)}\in\Fock$. Then $\Lambda'\in P^+$ and
$L(\Lambda')=\Usl s_{(0|\dots|0}$ is isomorphic to the integrable
highest weight $\Usl$-module of highest weight~$\Lambda'$. The bar
involution $\BarInvolution$ induces a unique $\Q$-linear
bar involution~$\BarInvolution$ on~$L(\Lambda')$ such that
$\overline{s_{\Lambda'}}=s_{\Lambda'}$ and
$\overline{u\cdot s_{\Lambda'}}=\overline{u}\cdot s_{\Lambda'}$, for all $u\in\Usl$.
Brundan and Kleshchev~\cite[Theorem~3.26]{BK:GradedDecomp} show that the
bar involution on~$L(\Lambda')$ extends to an involution
$\BarInvolution$ on~$\Fock$ and hence that the following holds.

\begin{Theorem}\label{T:Uglov}
   Suppose that $\bmu\in\Parts$. Then there is a unique bar-invariant
   element $G^\bmu\in\Fock$ such that
   \[ G^\bmu=\sum_{\blam\in\Parts} d_{\blam\bmu}(q)s_\blam,\]
   for some polynomials $d_{\blam\bmu}(q)\in\delta_{\blam\bmu}+q\Z[q]$.
\end{Theorem}

Using the coproduct $\Delta$ on $\Usl$, it is straightforward to show
that $\Fock\cong\Fock[-\kappa_\ell]\otimes\dots\otimes\Fock[-\kappa_1]$
as $\Usl$-modules. Uglov~\cite{Uglov} has proved a
more general version of these results where the Fock space does not
necessarily satisfy this tensor product decomposition.

Fayers~\cite{Fayers:LLT} used the tensor product decomposition
of~$\Fock$ to give an algorithm for computing the elements~$G^\bmu$
whenever $\bmu=(\mu^{(1)}|\dots|\mu^{(\ell)})$ where the partition
$\mu^{(r)}$ is $e$-restricted, for $1\le r\le \ell$. This set of
multipartitions contains $\Klesh'$. (A partition $\mu$ is
\textbf{$e$-restricted} if $\mu_k-\mu_{k+1}<e$, for $k\ge1$.)


Given $\blam\in\Parts$ define
$\beta_\blam=\sum_{A\in[\blam]}\alpha_{\res A}\in Q^+$.  The
combinatorial classification of the blocks of~$\Hn'$ and $\Sn'$ given
in~\cite{LM:AKblocks,Brundan:degenCentre} is equivalent to the statement
that the Specht modules, or Weyl modules, indexed by $\blam$ and $\bmu$
are in the same block if and only if $\beta_\blam=\beta_\bmu$.
Following \cite{BK:GradedKL} define the \textbf{defect} of~$\blam$ by
\begin{equation}\label{E:Defect}
   \defect\blam=(\Lambda',\beta_\blam)-\frac12(\beta_\blam,\beta_\blam)\in\N.
\end{equation}
By the remarks above, the defect is a block invariant. The defect is
easily seen to be equivalent to the combinatorial definition of the
\textbf{weight} $w(\blam)$ of~$\blam$ given by
\cite[\S2.1]{Fayers:AKweight}. We can now state Fayers'
conjecture~\cite{Fayers:LLT}.

\begin{Conjecture}[Fayers~\cite{Fayers:LLT}]\label{Con:Fayers}
    Let $\blam,\bmu\in\Parts$
    Then $\deg d_{\blam\bmu}(q)\le\defect\bmu$ and, moreover,
    $\deg d_{\blam\bmu}(q)=\defect\bmu$ only if $\bmu\in\Klesh$.
\end{Conjecture}

To prove this conjecture, we work in a graded setting where we can interpret
$d_{\blam\bmu}(v)$ as a graded decomposition number. Building on work of
Khovanov and Lauda~\cite{KhovLaud:diagI} and
Rouquier~\cite{Rouquier:QuiverHecke2Lie}, Brundan and
Kleshchev~\cite{BK:GradedKL} showed that (each block of) $\Hn$ is a
$\Z$-graded algebra. Extending this result, Stroppel and
Webster~\cite{StroppelWebster:QuiverSchur} and the
authors (when $e=\infty$)~\cite{HuMathas:QuiverSchurI} showed that the cyclotomic
Schur algebras admit a $\Z$-grading. Following Stroppel and
Webster~\cite{StroppelWebster:QuiverSchur}, we will work with the twisted
cyclotomic Schur algebra~$\Sn'$.

Let $\DSn$ and $\DHn$ be the basic graded algebras of $\Sn'$ and $\Hn'$,
respectively, and let $\GRMod\DSn$ and $\GRMod\DHn$ be the corresponding
categories of finite dimensional graded modules with homogeneous maps of
degree zero. (Unlike in the ungraded setting, graded basic algebras are
not uniquely determined up to isomorphism, as graded algebras. There is,
however, a unique grading on~$\DSn$ such that \autoref{T:Koszul} below
holds and by applying the graded Schur functor this fixes the grading
on~$\DHn$. More explicitly,~$\DSn$ is the graded endomorphism algebra of
$\bigoplus_{\blam\in\Parts}\DP^\blam$, where $\DP^\blam$ is the graded
projective cover of the corresponding self-dual graded simple module.)

The algebra $\DSn$ comes equipped with graded standard modules
$\DDelta^\blam$ and graded simple modules~$\DL^\blam$, where we fix the
grading on these modules by requiring that $\DL^\blam$ is graded
self-dual and that there is a homogeneous surjection
$\DDelta^\blam\surjection\DL^\blam$, for $\blam\in\Parts$.
Similarly, the graded algebra $\DHn$ has graded analogues of the Specht
modules $\DS^\blam$ and graded simple modules~$\DD^\bmu$ where we again
fix the gradings requiring that~$\DD^\bmu$ is graded self-dual and
that~$\DS^\bmu$ surjects onto $\DD^\bmu$, for $\blam\in\Parts$ and
$\bmu\in\Klesh$. In particular, this implies that $\DL^\bmu$ and
$\DD^\mu$ are both one dimensional modules concentrated in degree zero.

If $M=\bigoplus_{d\in\Z}M_d$ is a $\Z$-graded module let $\<\ \>$ be
the shift functor so that $M\<s\>$ is the $\Z$-graded module with
homogeneous component of degree $d$ being $(M\<s\>)_d=M_{d-s}$. Then the
\textbf{graded decomposition numbers} of $\DSn$ and $\DHn$ are
the Laurent polynomials
\[
  [\DDelta^\blam:\DL^\bmu]_q=\sum_{d\in\Z}[\DDelta^\blam:\DL^\bmu\<d\>]q^d
  \quad\text{and}\quad
  [\DS^\blam:\DD^\bmu]_q=\sum_{d\in\Z}[\DS^\blam:\DD^\bmu\<d\>]q^d,
\]
for $\blam,\bmu\in\Parts$. Abusing notation slightly, as in
\autoref{S:SchurPairs} there is a graded Schur functor
\[\Funn'\map{\GRMod\DSn}{\GRMod\DHn}\]
such that $\Funn'(\DDelta^\blam)\cong\DS^\blam$ and
$\Funn'(\DL^\bmu)\cong\DD^\bmu$ which, of course, is zero if
$\bmu\notin\Klesh'$.

\begin{Theorem}[\protect{%
     Stroppel-Webster~\cite[Theorem~7.11]{StroppelWebster:QuiverSchur}}]
  \label{T:StroppelWebster}
  Suppose that $K=\C$ and let $\blam,\bmu\in\Parts$. Then
\[[\DDelta^\blam:\DL^\bmu]_q=d_{\blam\bmu}(q).\]
Moreover, if $\bmu\in\Klesh$ then $[\DS^\blam:\DD^\bmu]_q=d_{\blam\bmu}(q)$.
In particular, $d_{\blam\bmu}(q)\in\delta_{\blam\bmu}+q\N[q]$, for all
$\blam,\bmu\in\Parts$.
\end{Theorem}

The graded decomposition numbers of $\DHn$ were first computed by
Brundan and Kleshchev~\cite{BK:GradedDecomp}. We need one
property of these polynomials, which was part of Fayers' motivation for
\autoref{Con:Fayers}.

\begin{Lemma}[\protect{%
  \cite[Remark~3.19]{BK:GradedDecomp} and \cite[Corollary~3.6.7]{Mathas:Singapore}}]
  \label{L:Defect}
  Suppose that $K=\C$, $\blam\in\Parts$ and $\bmu\in\Klesh'$. Then
  $0\le \deg d_{\blam\bmu}(q)\le\defect\bmu$. Moreover, the following
  are equivalent:
  \begin{enumerate}
    \item $\deg d_{\blam\bmu}(q)=\defect\bmu$,
    \item $d_{\blam\bmu}(q)=q^{\defect\bmu}$,
    \item $\blam=\bm^{-1}(\bmu)'$.
  \end{enumerate}
\end{Lemma}

Using properties of the polynomials $d_{\blam\bmu}(q)$ from
\autoref{T:Uglov}, Fayers~\cite[Corollary~2.4]{Fayers:AKweight} proved that if
$\bmu\in\Klesh'$ then there is a unique multipartition $\blam\in\Parts$
such that $d_{\blam\bmu}(q)=q^{w(\bmu)}$, where $w(\bmu)$ is his weight
function. Combining Fayers' result with \autoref{L:Defect} gives a
representation theoretic proof that $w(\bmu)=\defect\bmu$.

We need one more (deep) result from the graded representation theory
of~$\DSn$. We refer the reader to \cite[\S2]{BGS:Koszul} for the
definition of a \textit{Koszul} algebra.

\begin{Theorem}[\protect{Hu-Mathas~\cite{HuMathas:QuiverSchurI}
  ($e=\infty$) and Maksimau~\cite{Maksimau:QuiverSchur}}]\label{T:Koszul}
  Suppose that $K=\C$. Then the basic twisted cyclotomic Schur algebra $\DSn$
  is a Koszul algebra.
\end{Theorem}

The papers \cite{HuMathas:QuiverSchurI,Maksimau:QuiverSchur} actually
prove this result for the blocks of $\DSn$, however, this implies
\autoref{T:Koszul} because the direct sum of two Koszul algebras is again
Koszul.

Notice that \autoref{T:StroppelWebster} implies that $\DSn$ and $\DHn$ are
both positively graded algebras (strictly, non-negatively graded
algebras). Moreover, since $\DL^\bmu$ and $\DD^\bmu$ are both
concentrated in degree zero, this implies that each of the modules
$\DDelta^\blam$, $\DL^\blam$, $\DS^\blam$ and $\DD^\blam$ is positively
graded since $d_{\blam\bmu}(q)\in\delta_{\blam\bmu}+q\N[q]$. If
$M=\bigoplus_{d\ge0}M_d$ is a positively graded $\DSn$-module then its
\textbf{grading filtration} is
\[ M=\Gr_0M\supseteq \Gr_1 M\supset\dots\]
is given by $\Gr_r M=\bigoplus_{d\ge r}M_d$. Since $\DSn$ is a positively
graded algebra, $\Gr_r M$ is an $\DSn$-module for~$r\ge0$. In
particular, $\DDelta^\blam$ has a grading filtration.
For $\blam,\bmu\in\Parts$ write
\[  d_{\blam\bmu}(q) = \sum_{r\ge0}d_{\blam\bmu}^{(r)}q^r,
     \qquad\text{where }d_{\blam\bmu}^{(r)}\in\N.
\]
Let
$\DDelta^\blam=\rad^0\DDelta^\blam\supseteq\rad^1\DDelta^\blam
              \supseteq\cdots$
be the radical filtration of~$\DDelta^\blam$.

\begin{Corollary}\label{C:Radical}
  Suppose that $K=\C$ and let $\blam\in\Parts$. Then
  $\rad^r\DDelta^\blam=\Gr_r\DDelta^\blam$. Consequently,
  \[
          [\rad^r\DDelta^\blam/\rad^{r+1}\DDelta^\blam:\DL^\bmu]_q
          = d_{\blam\bmu}^{(r)}q^r,
          \qquad\text{ for all }r\ge0.
  \]
\end{Corollary}

\begin{proof}
  By \cite[Corollary~2.3.3]{BGS:Koszul}, any Koszul algebra is quadratic.
  Therefore, since $\DDelta^\blam/\rad\DDelta^\blam\cong\DL^\blam$ is
  simple and concentrated in degree zero, the radical
  filtration of~$\DDelta^\blam$ coincides with the grading filtration
  of~$\DDelta^\blam$ by \cite[Proposition~2.4.1]{BGS:Koszul}.
\end{proof}

We can now prove a stronger version of \autoref{Con:Fayers} (and of
\autoref{L:Defect}).

\begin{Theorem}[Fayers' Conjecture]\label{T:Fayers}
  Suppose $K=\C$ and that $\blam,\bmu\in\Parts$. Then
  $\deg d_{\blam\bmu}(q)\le\defect\bmu$ with equality only if $\bmu\in\Klesh'$.
  Moreover, the following are equivalent:
  \begin{enumerate}
     \item $\deg d_{\blam\bmu}(q)=\defect\bmu$,
     \item $d_{\blam\bmu}(q)=q^{\defect\bmu}$,
     \item $\blam=\bm^{-1}(\bmu)'$,
     \item $\soc\DDelta^\blam=\DL^\bmu\<\defect\bmu\>$,
     \item $\soc\DS^\blam=\DD^\bmu\<\defect\bmu\>$.
  \end{enumerate}
\end{Theorem}

\begin{proof}
  First observe that by the comments before \autoref{E:Defect}, if
  $d_{\blam\bmu}(q)\ne0$ then $\defect\blam=\defect\bmu$. Fix
  $\bmu\in\Parts$ such that $d_{\blam\bmu}(q)\ne0$ and $d=\deg d_{\blam\bmu}(q)$
  is maximal in the sense that $\deg d_{\blam\bnu}(q)\le d$ for all
  $\bnu\in\Parts$. By the maximality of~$d$, $\Gr_{d+1}\DDelta^\blam=0$
  and $\DL^\bmu\<d\>$ is a summand of
  $\Gr_d\DDelta^\blam=\rad^d\DDelta^\blam$, where the last equality
  comes from \autoref{C:Radical}. Consequently, $\DL^\bmu\<d\>$ is a
  summand of the socle of~$\DDelta^\blam$. Forgetting the gradings, this
  implies that $L_\prime^\bmu$ is contained in socle of
  $W^\blam_\prime$.  Therefore, $\bmu\in\Klesh'$ by
  \autoref{T:DeltaSocle}. As~$\bmu\in\Klesh'$ and
  $\DL^\bmu\<d\>\subseteq\soc\DDelta^\blam$ whenever
  $d=\deg d_{\blam\bmu}(q)$ is maximal the theorem now follows by
  applying \autoref{L:Defect}.
\end{proof}

Notice that by interchanging the roles of $\charge$ and $\charge'$,
\autoref{T:Fayers} becomes a result about~$\Sn$-modules. By
\autoref{T:Fayers} if $\blam'\notin\Klesh$ then the socle of $\DDelta^\blam$
will not have a component in degree~$\defect\blam$ and, \textit{a
priori}, $\soc\DDelta^\blam$ is not necessarily homogeneous. On the
other hand, by \autoref{T:DeltaSocle} and \autoref{T:Fayers}, if $d\in\Z$ and~$\DL^\bmu\<d\>$ appears
in~$\soc\DDelta^\blam$ then $0\le d\le\defect\blam$ and~$\bmu\in\Klesh'$.

\section{Homomorphisms between Weyl modules and Specht modules}
  \label{S:Homomorphisms}

  In our final section we return to ungraded representation theory of
  the cyclotomic Schur algebras~$\Sn$ and we prove some results relating
  the hom-spaces between Specht modules and Weyl modules. The main
  result in this section is a cyclotomic analogue of the classical
  Carter-Lusztig Theorem~\cite[Theorem~3.7]{CarterLusztig}. The
  $q$-analogue of this result in level one was proved by Dipper and
  James~\cite{DJ:qWeyl}.

\begin{Lemma} \label{L:injection} Suppose that  $X, Y\in\rDelMod\Sn$.
  Then the Schur functor~$\Funn$ induces an injection
  \[
      \Hom_{\Sn}(X,Y)\injection\Hom_{\Hn}(\Funn(X),\Funn(Y));
                 \phi\mapsto\Funn(\phi).
  \]
\end{Lemma}

\begin{proof} Since $\Funn$ is exact it suffices to consider the case when
  $X=W^{\blam}$ and $Y=W^{\bmu}$, for some $\blam,\bmu\in\Parts$.
Suppose that $f\in\Hom_{\Sn}(W^{\blam},W^{\bmu})$ and that $\Funn(f)=0$.
It follows that
\[
\Funn(\im(f))=\im(\Funn(f))=0.
\]
Now, if $f\ne0$ then $\im f\ne0$ so there exists $\bmu\in\Klesh$ such
that $[\im f:L^\bmu]\ne0$ by \autoref{T:DeltaSocle}. Therefore,
$\Funn(\im f)\ne0$, showing that $\Funn(f)=0$ if and only if $f=0$.
\end{proof}

The main result of this section gives sufficient conditions for the map
of \autoref{L:injection} to be an isomorphism of vector spaces. In the
language of Rouquier~\cite[Definition~4.37]{Rouquier:Schur}, our result
gives sufficient conditions for $\Sn$ to be a \textit{$0$-faithful cover}
of~$\Hn$.

By \autoref{L:injection} there is an injection
$\Hom_{\Sn}(W^\blam,W^\bmu)\injection\Hom_{\Hn}(S^\blam,S^\bmu)$, for
all $\blam,\bmu\in\Parts$. We want to recast this in the framework
developed in \autoref{S:SchurPairs} to prove \autoref{L:Annihilators}. To this
end, for $\blam\in\Parts$ let $\plam=\varphi_{\Tlam\Tlam}$ be
the identity map on $M^\blam$ and define $\ELam=\Sn\pom$, where
$\omega=(0|\dots|0|1^n)\in\Parts$.

Using \autoref{C:DualSpechts} to take duals
gives an isomorphism
$\Hom_{\Hn}(S^\blam,S^\bmu)\cong\Hom_{\Hn}(S_\bmu,S_\blam)$.
We work with homomorphisms between dual Specht modules because
this better fits the framework developed in \autoref{S:SchurPairs}.

Fix $\bmu\in\Parts$ and recall from \autoref{S:HeckeAlgebras} that
$z_\bmu=m_\bmu T_{w_\bmu} n_\bmu$. Let $\pmu$ be the identity map
on~$M^\bmu=m_\bmu\Hn$ and let $\pimu\map\Hn M^\mu$ be the natural
surjection given by $\pi_\bmu(h)=m_\bmu h$ for $h\in\Hn$. In fact,
$\pmu=\varphi_{\bT^\bmu\bT^\bmu}$ and $\pimu=\varphi_{\bT^\bmu\t^\bmu}$
are both elements of~$\Sn$. The isomorphism $\Hn\cong\End_{\Hn}(\Hn)$,
which maps $h\in\Hn$ to left multiplication by~$h$, identifies $\Hn$
with the subalgebra $\pom\Sn\pom$, where
$\pom=\varphi_{\bT^\bom\bT^\bom}$ is the identity map on~$\Hn$. Identify
$\Hn$ with its image under this map.  Let $\zmu = \pimu
T_{w_\bmu}n_\bmu$. Then $\zmu\in\pmu\Sn\pom$. Define $\Delta^\bmu =
\Sn\zeta_\bmu$.

\begin{Lemma}\label{L:LR_S}
  Suppose that $\zeta\in\Sn\pom$. Then
  $\LR_{\Sn}(\zeta)=\Sn\pom\cap\L_{\Sn}\bigl(\R_{\Hn}(\zeta)\bigr)$.
\end{Lemma}

\begin{proof}
    First observe that $\zeta(1-\pom)=0$ because $\zeta\in\Sn\pom$, so if
    $s\in\LR_{\Sn}(\zeta)$ then $s(1-\pom)=0$. Hence,
    $\LR_{\Sn}(\zeta)\subseteq\Sn\pom$ and, consequently,
    $\LR_{\Sn}(\zeta)\subseteq\Sn\pom\cap\L_{\Sn}\bigl(R_{\Hn}(\zeta)\bigr)$.
    To prove the reverse inclusion, let
    $x\in\Sn\pom\cap\L_{\Sn}\bigl(\R_{\Hn}(\zeta)\bigr)$ and fix
    $s\in\R_{\Sn}(\zeta)$. Then
    $0=\zeta s = \zeta\pom s = \zeta\pom s\pilam$,
    for all $\blam\in\Parts$ (of course, $\pom
    s\pilam$ could be zero). By definition, $\pom s\pilam\in\Hn$, so
    $x\pom s\pilam=0$ since
    $x\in\L_{\Sn}\bigl(\R_{\Hn}(\zeta)\bigr)$. Hence, $x\pom s\plam=0$ since
    $\pilam=\ilam\pilam$ is surjective, for $\blam\in\Parts$, so
    \[xs = x\pom s = \sum_{\blam\in\Parts}x\pom s\plam =0.\]
    Hence, $x\in\LR_{\Sn}(\zeta)$ as we needed to show.
\end{proof}

As in \autoref{S:SchurPairs}, the $*$-isomorphism of $\Hn$ induces an
anti-isomorphism of~$\Sn$, which we also call~$*$ (see after
\autoref{D:SelfDualSchurPair}). Like we did in \autoref{D:rightRdfn}, if $X$ is a right $\Sn$-module let $X_L$
be the left $\Sn$-module that is equal to $X$ as a vector space and with
left action given by $s\cdot x = xs^*$, for $s\in\Sn$ and $x\in X$.
Similarly, if $X$ is a right $\Hn$-module then let $X_L$ be the
corresponding left $\Hn$-module.

\begin{Lemma}\label{L:CarterLusztig}
    Suppose that $\blam,\bmu\in\Parts$. Then
    $S_\blam\cong\zeta_\blam\Hn$, as right $\Hn$-modules,
    $\Delta^\blam\cong W^\blam_L$, as left $\Sn$-modules, and
    there are vector space isomorphisms,
    \begin{align*}
      \Hom_{\Sn}(W^\blam,W^\bmu)  &\cong \Hom_{\Sn}(\Dlam,\Dlam[\bmu])
      \cong \Dlam[\bmu]\cap\RL_{\Sn}(\zlam)
    \intertext{and}
      \Hom_{\Hn}(S^\blam,S^\bmu) &\cong \Hom_{\Hn}(S_\bmu,S_\blam)
                                  \cong\LR_{\Sn}(\zmu)\cap S_\blam.
    \end{align*}
\end{Lemma}

\begin{proof}
  The map $m_\blam h\mapsto \pilam h$, for $h\in\Hn$, defines a
  $\Hn$-module isomorphism $M^\blam\cong\varphi_\blam\Sn\varphi_\bom$
  that sends~$z_\blam$ to $\zlam$. Therefore,
  $S_\blam\cong\zeta_\blam\Hn$ by
  \autoref{L:SpechtSubmodules}. Similarly, $\Dlam\cong(W^\blam)^*$ as
  left $\Sn$-modules by
  \autoref{P:WeylSubmodules}. Now consider the two hom-spaces. First,
  since $\Dlam[\bnu]\cong W^\bnu_L$ for $\bnu\in\Parts$,
  \[
  \Hom_{\Sn}(W^\blam,W^\bmu)\cong\Hom_{\Sn}(\Dlam,\Dlam[\bmu])
       \cong\Hom_{\Sn}(\Sn\zlam,\Dlam[\bmu])
       \cong\RL_{\Sn}(\zlam)\cap\Dlam[\bmu],
  \]
  using \autoref{L:Annihilators}(a). Finally,
  if $\bnu\in\Parts$ then $(S^\bnu)^*\cong S_\bnu$ so taking duals,
  \begin{align*}
       \Hom_{\Sn}(S^\blam,S^\bmu)
           &\cong \Hom_{\Sn}(S_\bmu,S_\blam)
            \cong \Hom_{\Sn}(\zmu\Hn,\zlam\Hn).
  \end{align*}
  To complete the proof it is enough to show that
  $\Hom_{\Hn}(\zmu\Hn,\zlam\Hn)\cong\LR_{\Sn}(\zmu)\cap S_\blam$.
  By \autoref{L:LR_S}, if $f\in\Hom_{\Hn}(\zmu\Hn,\zlam\Hn)$ then
  $f(\zmu)\in\LR_{\Sn}(\zmu)$. Conversely, if $z\in\LR_{\Sn}$ then
  \autoref{L:LR_S} implies that there is a well-defined $\Hn$-module
  homomorphism $f\map{\zmu\Hn}\zlam\Hn$ given by $f(\zmu h)=zh$, for
  $h\in\Hn$. Hence, $\Hom_{\Hn}(\zmu\Hn,\zlam\Hn)\cong\LR_{\Sn}(\zmu)\cap S_\blam$ via the map $f\mapsto f(\zmu)$.
\end{proof}

So, to prove that
$\Hom_{\Sn}(W^\blam,W^\bmu)\cong\Hom_{\Sn}(S^\blam,S^\bmu)$ it is enough
to show that
\[\Dlam[\bmu]\cap\RL_{\Sn}(\zlam)=\LR_{\Hn}(\zmu)\cap\zlam\Hn.\]
To compare these hom-spaces we need to describe $\RL_{\Sn}(\zlam)$ and
$\LR_{\Sn}(\zmu)$.  To compute $\LR_{\Sn}(\zmu)$ we prove a
series of results about the ``left Weyl modules'' $\Dlam$, for
$\blam\in\Parts$. By \autoref{L:CarterLusztig}, $W^\blam=\Dlam_R$ so all
of these results translate into statements about Weyl modules, which we
leave as a exercise for the reader.

First, we need a fact that is of the folklore for~$\Hn$ but it does not seem
to be in the literature. In level one, this is a result of Dipper and
James~\cite[Lemma~4.1]{DJ:reps}.

\begin{Lemma}\label{L:mnideals}
  Suppose that $m_\bsig\Hn n_\bnu\ne0$, for $\bsig,\bnu\in\Parts$. Then
  $\bnu\gedom\bsig$.
\end{Lemma}

\begin{proof}
  By \cite[Theorem~4.14]{DJM:cyc}, and as discussed in \autoref{S:SchurAlgebras},
  $M^\bsig=m_\bsig\Hn$ has basis
  $\set{m_{\bS\t}|\bS\in\RSStd(\btau,\bsig), \t\in\Std(\btau)
             \text{ for }\btau\in\Parts}$.
  Therefore, $m_\bsig\Hn n_\bnu\ne0$ only if there exist tableaux
  $\bS\in\RSStd(\btau,\bnu)$ and $\t\in\Std(\btau)$, for some $\btau\in\Parts$,
  such that $m_{\bS\t}n_\bnu\ne0$. Hence, $\t_\bnu\gedom\t$ by
  \autoref{P:mnDual}. Consequently, $\t_\bnu\gedom\t_\btau$, so that
  $\bnu\gedom\btau\gedom\bsig$, as required.
\end{proof}

Recall the element $n_\bmu\in\Hn$ from \autoref{S:HeckeAlgebras}. Let
$\theta_\bmu$ be the natural embedding of~$N^\bmu=n_\bmu\Hn$ into
$\Hn=M^{\bom}$ and set $E^\bmu=\Sn\theta_\bmu$. Then $E^\bmu$
is a left $\Sn$-submodule of $\Hom_{\Hn}(N^\bmu,M)$.  (The right $\Sn$-module
$\bigoplus_\bmu E^\bmu_R$ is the full tilting module underpinning the
proof of \autoref{T:RingelDuality}; see \cite[Theorem~6.18]{M:tilting}.)

The remaining results in this section depend on \cite{M:tilting}, which
only considers the non-degenerate Hecke algebras, so \textit{henceforth
we assume that $\xi^2\ne1$}. In view of \autoref{D:Hecke} and
\cite[Corollary~2.10]{HuMathas:GradedInduction}, the arguments of
\cite{M:tilting} should extend to the degenerate case, however, these
results do not appear in the literature so we cannot use them.

For the proof of the next result recall from \autoref{S:SchurAlgebras} that
$\RSStd(\bnu,\bmu)$ and $\CSStd(\bnu,\bmu)$, respectively, are the sets
of row and column semistandard $\bnu$-tableau of type~$\bmu$. If
$\bT\in\CSStd(\bnu,\bmu)$ let $\dot\bT$ be the unique standard
$\bnu$-tableau such that $\col_\bmu(\dot\bT)=\bT$ and
$\ell(d(\dot\bT))\le\ell(d(\t))$ whenever $\t$ is standard and
$\col_\bmu(\t)=\bT$. Mirroring the definition of the $\varphi$-basis of~$\Sn$,
for $\bU\in\CSStd(\btau,\bmu)$ and $\bV\in\CSStd(\btau,\bnu)$, where
$\btau\in\Parts$, define $\psi_{\bU\bV}\in\Hom_{\Hn}(N^\bnu,N^\bmu)$ by
$\psi_{\bU\bV}(n_\bnu h)=n_{\bU\bV}h$, for $h\in\Hn$ and
$n_{\bU\bV}=\sum_{\u,t}(-\xi)^{-\ell(d(\u))-\ell(d(\v))}n_{\u\v}$, where
the sum is over all standard tableaux such that $\col_\bmu(\u)=\bU$ and
$\col_\bnu(\v)=\bV$. Then \autoref{T:CellularSAlgebras} and
\autoref{P:Mullineux} imply that $\set{\psi_{\bU\bV}}$ is a cellular basis
of $\Sn(N)$.

The next result describes the left Weyl module $\Dlam[\bmu]$ as a
kernel intersection.

\begin{Lemma}\label{L:Kernel}
  Suppose that $\xi^2\ne1$ and that $\bmu\in\Parts$ and for $\bnu\in\Parts$ let
  $\Emunu=\Hom_{\Sn}(E^\bmu,E^\bnu)$. Then
  \[ \Dlam[\bmu]=\bigcap_{\bnu\not\gedom\bmu}
        \bigcap_{\Psi\in\Emunu}\ker\Psi.\]
\end{Lemma}

\begin{proof}
   In view of \cite[Proposition~6.4]{M:tilting}, a basis of $E^\bmu$ is given by
   the maps
   \[\set{\theta_{\bS\bT}|\bS\in\RSStd(\bsig,\bal), \bT\in\CSStd(\bsig,\bmu)
         \text{ for }\bal,\bsig\in\Parts}\]
   where $\theta_{\bS\bT}(n_\bmu h) = m_{\bS\dot\bT}n_{\bmu}h$, for
   $h\in\Hn$. As noted in \cite[Theorem~6.5]{M:tilting}, this implies
   that $E^\bmu$ has a filtration by left Weyl modules $\Dlam[\bsig]$ such that
   $(E^\bmu: \Dlam[\bsig])=\#\CSStd(\bsig,\bmu)$. In particular,
   $(E^\bmu: \Dlam[\bsig])\ne0$ only if $\bmu\gedom\bsig$. By
   \cite[Proposition~7.1]{M:tilting}, the hom-space $\Emunu$ has basis
   \[\set{\Psi_{\bU\bV}|\bU\in\CSStd(\btau,\bmu)\text{ and }
           \bV\in\CSStd(\btau,\bnu) \text{ for }\btau\in\Parts},
   \]
   where $\Psi_{\bU\bV}(\theta)=\theta\circ\psi_{\bU\bV}$ for
   $\theta\in E^\bmu$. (As in \autoref{R:nbasislabels}, the notation used
   in \cite{M:tilting} is slightly different with what we use here in that we are
   working with left $\Hn$-modules here, so some care
   must be taken when comparing our results with those
   of~\cite{M:tilting}.) Armed with these facts we can prove the lemma.

   Fix $\bS\in\RSStd(\bsig,\bal)$, $\bT\in\CSStd(\bsig,\bmu)$,
   $\bU\in\CSStd(\btau,\bmu)$ and $\bV\in\CSStd(\btau,\bnu)$. Then the
   $\Hn$-module homomorphism $\Psi_{\bU\bV}(\theta_{\bS\bT})\in E^\bnu$
   is completely determined by
   \begin{equation}\label{E:PsiTheta}
      \Psi_{\bU\bV}(\theta_{\bS\bT})(n_\bnu)
          =(\theta_{\bS\bT}\circ\psi_{\bU\bV})(n_\bnu)
          =\theta_{\bS\bT}(n_{\bU\bV})
          =\theta_{\bS\bT}(n_\bmu h_{\bU\bV})
          =m_{\bS\dot\bT}n_\bmu h_{\bU\bV}
          =m_{\bS\dot\bT}n_{\bU\bV},
   \end{equation}
   where $n_{\bU\bV}=n_\bnu h_{\bU\bV}$.
   In particular, $\Psi_{\bU\bV}(\theta_{\bS\bT})\ne0$ only if
  $\bnu\gedom\bsig$ by \autoref{L:mnideals}.

   By construction, the module $\Dlam[\bmu]$ has basis
   $\set{\theta_{\bS\bT_\bmu}|\bS\in\RSStd(\bmu,\bal)
                 \text{ for }\bal\in\Parts}$,
   so $\Dlam[\bmu]$ is a submodule of $E^\bmu$. Indeed,
   $\zlam[\bnu]=\theta_{\bT^\bmu\bT_\bmu}$ so the map
   $\phi\mapsto\phi\theta_{\bT^\bmu\bT_\bnu}$ is an injection. Hence,
   taking $\bsig=\bmu$ and $\bT=\bT_\bmu$ in~\autoref{E:PsiTheta},
   $\Dlam[\bmu]\subseteq\bigcap_{\bnu\not\gedom\bmu}
      \bigcap_{\Psi\in\Emunu}\ker\Psi$.

   Conversely, if $\bsig\ne\bmu$ and there exist tableaux
   $\bS\in\RSStd(\bsig,\bal)$ and $\bT\in\CSStd(\bsig,\bmu)$ then
   $\theta_{\bS\bT}\in E^\bmu$ and~$\bmu\gdom\bsig$
   since~$\CSStd(\bsig,\bmu)\ne\emptyset$. Moreover,
   $\Psi_{\bT\bT_\bsig}(\theta_{\bS\bT})\ne0$ by \autoref{E:PsiTheta}
   and \autoref{P:mnDual}. That is,
   $\theta_{\bS\bT}\notin\ker\Psi_{\bT\bT_\bsig}$. This
   calculation is independent of~$\bS$ and, in fact, it shows that
   $\Psi_{\bT\bT_\bsig}\bigl(\sum_\bS c_\bS\theta_{\bS\bT}\bigr)\ne0$
   for any scalars $c_\bS\in K$, which are not all zero. More generally,
   suppose that $\theta=\sum_{\bU,\bV}c_{\bU\bV}\theta_{\bU\bV}\in
   E^\bmu$, for some scalars $c_{\bU\bV}\in K$ such that
   $c_{\bU\bV}\ne0$ only if $\Shape(\bV)\gdom\mu$. Fix $(\bS,\bT)$ such
   that $c_{\bS\bT}\ne0$ and where $\bT$ is minimal in the sense that
   $\bT\ntriangleright\bV$ whenever $c_{\bU\bV}\ne0$ for some
   $(\bU,\bV)$. Combining \autoref{P:mnDual} with what we
   have just shown, it follows that
   $\Psi_{\bT\bT_\bsig}(\theta)
      =\sum_{\bS\in\RSStd(\bsig)}c_{\bS\bT}\Psi_{\bT\bT_\bsig}(\theta_{\bS\bT})
       \ne0$.

   Combining the last two paragraphs completes the proof.
\end{proof}

\begin{Corollary}\label{C:WeylIntersection}
    Suppose that $\xi^2\ne1$ and that $\bmu\in\Parts$. Then
    $\Dlam[\bmu]=E^\bmu\cap\LR_{\Sn}(\zmu)$.
\end{Corollary}

\begin{proof}
  We have $\Dlam[\bmu]=\Sn\zmu\subseteq E^\bmu\cap\LR_{\Sn}(\zmu)$ since
  $\zmu\in E^\bmu$. On the other hand, if
  $\bnu\not\gedom\bmu$ then
  every homomorphism $\Psi\in\Emunu$ is given by right multiplication
  by some element~$\psi$ by \autoref{E:PsiTheta}. Moreover, $\Dlam[\bmu]\subseteq\ker\Psi$
  \autoref{L:Kernel}, so if $\psi\in\Sn$ then $\zmu\psi=\Psi(\zmu)=0$. Therefore,
  if $x\in E^\bmu\cap\LR_{\Sn}(\zmu)$ then $\Psi(x)=x\psi=0$. Therefore,
  $E^\bmu\cap\LR_{\Sn}(\zmu)\subseteq\ker\Psi$, for all
  $\Psi\in\Emunu$. The corollary now follows by \autoref{L:Kernel}.
\end{proof}

\begin{Lemma}\label{L:nnuAnnihilator}
  Suppose that $\xi^4\ne1$ and $\kappa_r\not\equiv\kappa_s\pmod{e\Z}$,
  for $1\le r<s\le\ell$. Then $\LR_{\Sn}(n_\bmu)=E^{\bmu}$.
\end{Lemma}

\begin{proof} On the one hand, we have that
  $E^{\bmu}\subseteq \LR_{\Sn}(\theta_\bmu)=\LR_{\Sn}(n_\bmu)$ by definition.  On the other hand, under exactly
  these assumptions, \cite[Corollary~6.11]{M:tilting} says that
  \begin{align*}
    E^\bmu  & =\Hom_{\Hn}(M,N^\bmu)\cong\Hom_{\Hn}(N_L^{\bmu},M_L)
              \cong\Hom_{\Hn}\bigl(\pom(\Sn\theta_\bmu),\pom\Sn\bigr)\\
            &\cong \Hom_{\Hn}\bigl(\pom\Sn\otimes_{\Sn}\Sn\theta_\bmu, \pom\Sn\bigr)
             \cong \Hom_{\Sn}\bigl(\Sn\theta_\bmu,\Hom_{\Hn}(\pom\Sn, \pom\Sn)\bigr)\\
            &\cong\Hom_{\Sn}(\Sn\theta_\bmu,\Sn).
   \end{align*}
  Therefore, $E^\bmu=\LR_{\Sn}(\theta_\bmu)\cap\Sn=\LR_{\Sn}(n_\bmu)\cap\Sn=\LR_{\Sn}(n_\bmu)$ by
  \autoref{L:Annihilators}(a).
\end{proof}

\begin{Lemma}\label{L:KleshAnnihiator}
  Suppose that $\xi^2\ne1$ and $\bmu'\in\Klesh'$. Then $\LR_{\Sn}(\zmu)\subseteq E^{\bmu}$.
\end{Lemma}

\begin{proof}
  By \autoref{C:SimpleIsomorphism}, $D_{\bmu}\neq 0$ and so $n_{\bmu}\Hn
  n_{\bmu}\notin\Hn^{\lhd\bmu}$.  In particular, $z_\bmu\Hn n_\bmu\ne0$,
  so $\zmu\Hn n_\bmu\ne0$. So, there exists $h\in\Hn$ and $0\ne c\in K$
  such that $\zmu hn_\bmu=c \zmu$, or equivalently, $\zmu(c-hn_\bmu)=0$.
  Suppose that $\phi\in\LR_{\Sn}(\zmu)$. Then $\phi(c-hn_\bmu)=0$.
  Therefore, $c\phi=\phi hn_\bmu\in\Sn n_\bmu$, so that $\phi\in \Sn n_\bmu$
  since~$c\ne0$.  Hence, by \autoref{L:nnuAnnihilator}, $\LR_{\Sn}(\zmu)\subseteq \Sn n_\bmu\subseteq \LR_{\Sn}(n_\bmu)=E^{\bmu}$ as
  required.
\end{proof}

Before we can prove the main result of this section we need some
analogous results for $\RL_{\Sn}(\zlam)$. We start with an analogue of
\autoref{L:Kernel} for $S_\blam$.

\begin{Lemma}\label{L:SKernelIntersection}
  Suppose that $\xi^2\ne1$ and that $\blam\in\Parts$. Then
  \[ \zlam\Hn = \set{m\in\pilam\Hn|\phi m =0 \text{ for all } \phi\in\plam[\bnu]\Sn\plam
                  \text{ for }\bnu\in\Parts\text{ with }\blam\not\gedom\bnu}.\]
\end{Lemma}

\begin{proof}
  The argument is similar to the proof of \autoref{L:Kernel} so we just
  sketch the proof. Let
  \[ X_\blam= \set{m\in\pilam\Hn|\phi m =0 \text{ for all } \phi\in\plam[\bnu]\Sn\plam
                  \text{ for }\bnu\in\Parts\text{ with }\blam\not\gedom\bnu},
  \]
  so we need to show that $\zlam\Hn=X_\blam$.  In view of
  \cite[Proposition~5.9]{M:tilting}, $\pilam\Hn$ has basis\footnote{Note that the element $n_{\s\t}$ in this paper corresponds to the element $n_{\s'\t'}$ in the notation of \cite{M:tilting}. This accounts for the difference between our description of this basis and \cite[Proposition~5.9]{M:tilting}.}
  \[
    \set{\pilam n_{\ubs\t}|\bS\in\RSStd(\bmu,\blam)\text{ and }\t\in\Std(\bmu)
          \text{ for some }\bmu\in\Parts}
  \]
  and $\zlam\Hn$ is the $\Hn$-submodule of $\pilam\Hn$ with basis
  $\set{\zlam T_{d(\t)}|\t\in\Std(\blam)}$, where $\ubs$ is the minimal standard $\blam$-tableau (under the dominance $\lhd$) such that
  $\row_{\bmu}(\ubs)=\bmu$. Observe that
  $\zlam T_{d(\t)}=\pilam n_{\bT^{\blam}\t}$, for $\t\in\Std(\blam)$. If $\bnu\in\Parts$ and
  $\bU\in\RSStd(\brho,\bnu), \bV\in\RSStd(\brho,\blam)$ and $\t\in\Std(\bmu)$ then
  \[
  \phiUV\zlam T_{d(\t)} (1) \in m_\bnu \Hn n_\blam T_{d(\t)}.
  \]
  Consequently, if $\blam\not\gedom\bnu$ then $\phiUV\zlam T_{d(\t)}=0$
  by \autoref{L:mnideals}. Therefore, $\zlam\Hn\subseteq X_\blam$. To
  prove the reverse inclusion, if $\bmu\ne\blam$ and $\bS\in\RSStd(\bmu,\blam)$ then
  $\bmu\gdom\blam$ and $\varphi_{\bT^\bmu\bS}\pilam
  n_{\ubs\t}(1)\ne0$ by \autoref{P:mnDual}. Arguing as in the
  last paragraph of \autoref{L:Kernel} now completes the proof.
\end{proof}

\begin{Proposition}\label{P:RLzlam}
  Let $\blam\in\Parts$ and suppose that $\xi^2\ne1$. Then $\zlam\Hn = \RL_{\Sn}(\zlam)\cap\Sn\pom$.
\end{Proposition}

\begin{proof}
  Certainly, $\zlam\Hn=\zlam\Sn\cap\Sn\pom\subseteq\RL_{\Sn}(\zlam)\cap\Sn\pom$.
  Conversely, suppose that $x\in\RL_{\Sn}(\zlam)\cap\Sn\pom$. Then
  $x\in\plam\Sn\pom=\plam\Hn$ since $(1-\plam)\zlam=0$. If
  $\blam\not\gedom\bnu$ and $\bS\in\RSStd(\brho,\bnu), \bT\in\RSStd(\brho,\blam)$ then
  $\phiST\zlam\in m_\bnu\Hn n_\blam$ so that
  $\phiST\zlam=0$ since $m_\bnu\Hn n_\blam=0$ by \autoref{L:mnideals}.
  Therefore, $\phiST x=0$ so that $x\in\zlam\Hn$ by
  \autoref{L:SKernelIntersection}. The result follows.
\end{proof}

Finally we can prove our cyclotomic generalisation of the
Carter-Lusztig Theorem~\cite{CarterLusztig}.

\begin{Theorem}\label{T:CarterLusztig}
  Let $\blam,\bmu\in\Parts$ and assume that $\xi^2\ne1$ and that either
  \begin{enumerate}
    \item $\bmu'\in\Klesh'$, or
    \item $\xi^2\ne-1$ and $\kappa_r\not\equiv\kappa_s\pmod{e\Z}$,
          for $1\le r<s\le\ell$.
  \end{enumerate}
  Then
  $
      \Hom_{\Sn}(W^\blam,W^\bmu)\cong\Dlam[\bmu]\cap \zlam\Hn
                        \cong \Hom_{\Hn}(S^\blam,S^\bmu)
  $
  as vector spaces.
\end{Theorem}

\begin{proof}
  By \autoref{L:CarterLusztig} and \autoref{P:RLzlam},
  $\Hom_{\Sn}(W^\blam,W^\bmu)\cong\Delta^\bmu\cap \zlam\Hn$ and
  $\Hom_{\Hn}(S^\blam,S^\bmu)\cong\L_{\Sn}(\R_{\Hn}(\zmu))\cap \zlam\Hn$, so it is
  enough to show that $\LR_{\Sn}(\zmu)=\Dlam[\bmu]$, under the assumptions of
  the theorem. By \autoref{C:WeylIntersection},
  $\Dlam[\bmu]=E^\bmu\cap\LR_{\Sn}(\zmu)$, so it is enough to
  show that $\LR_{\Sn}(\zmu)\subseteq E^\bmu$ by \autoref{L:LR_S}. If
  $\bmu'\in\Klesh'$ then this is immediate from
  \autoref{L:KleshAnnihiator}. On the other hand, if the conditions
  in~(b) hold then $E^\bmu=\LR_{\Sn}(n_\bmu)$ by
  \autoref{L:nnuAnnihilator}, but
  $\LR_{\Sn}(\zmu)\subseteq\LR_{\Sn}(n_\bmu)=E^\bmu$, since
  $\zmu=\pilam[\bmu] T_{w_\bmu}n_\bmu$, so again the result follows.
\end{proof}

By \autoref{L:injection}, the Schur functor induces an injective map
$\Hom_{\Sn}(W^\blam,W^\bmu)\injection\Hom_{\Hn}(S^\blam,S^\bmu)$.
Under the assumptions of \autoref{T:CarterLusztig} this map is an isomorphism.
If neither of the conditions in~(a) and~(b) in \autoref{T:CarterLusztig} hold
then it is not difficult to find examples where these two hom-spaces are
not isomorphic. Examples that show that the conclusions of
\autoref{T:CarterLusztig} do not hold in general are easy to construct
starting from easy observations that if $\xi^2=-1$ then
$S^{(2)}\cong S^{(1^2)}$ and if $\kappa_r=\kappa_s$ then
$S^{\boldsymbol\eta_r}\cong S^{\boldsymbol\eta_s}$, where
$1\le r<s\le\ell$ and
$\boldsymbol\eta_t=(\eta_t^{(1)},\dots,\eta_t^{(\ell)})$ is the multipartition with
$\eta_t^{(k)}=(1)$ if $t=k$ and $\eta_t^{(k)}=(0)$ otherwise. Compare
with \cite[Remark~6.10]{M:tilting}.

The assumption that $\xi^2\ne1$ in \autoref{T:CarterLusztig} is almost
certainly unnecessary. As stated above, we include it because
\cite{M:tilting} does not consider the degenerate Hecke algebras,
although as far as we have checked the arguments from~\cite{M:tilting} also apply when
$\xi^2=1$.  Under the assumptions in part~(b) the isomorphism of
\autoref{T:CarterLusztig} is implied by \cite[Theorem
6.6]{Rouquier:Schur}.  As here, Rouquier's argument relies on results
from~\cite{M:tilting} and, in fact, the conditions in
\autoref{T:CarterLusztig}(b) come from \cite[Theorem~6.9]{M:tilting}.
We state the full result here because almost no extra effort is required
to prove \autoref{T:CarterLusztig} under the assumptions in parts~(a)
and~(b).


\end{document}